\documentclass[11pt]{article}
\pdfoutput=0
\usepackage{amsfonts, latexsym, amssymb,amsmath,amstext}
\usepackage[dvips]{graphics,color}
\usepackage[mathscr]{eucal}
\newenvironment{proof}{\par\noindent{\bf Proof.}\ }{\hfill$\Box$\par\medskip}
\newtheorem{theorem}{Theorem}[section]

\newcommand{\equal}{&\!\!\!=\!\!\!&}

\newcommand{\CC}{\mathbb C}
\newcommand{\PP}{\mathbb P}
\newcommand{\QQ}{\mathbb Q}
\newcommand{\fK}{{\Bbb F}}
\newcommand{\hpg}[5]{{}_{#1}\mbox{\rm F}_{\!#2}\!
  \left(\left.{#3 \atop #4}\right| #5 \right) }
\makeatletter
\@addtoreset{equation}{section}

\makeatother
\newcommand{\whilb}[2]{$\textstyle{#1}\atop#2$}

\title{Computation of highly ramified coverings} % for algebraic Painlev\'e VI functions}
\author{Raimundas Vid\= unas\thanks{Supported by the 21 Century COE Programme
"Development of Dynamic Mathematics with High Functionality" of the Ministry
of Education, Culture, Sports, Science and Technology of Japan. E-mail:
vidunas@math.kyushu-u.ac.jp}\; and\; Alexander~V.~Kitaev\thanks{Supported by
JSPS grant-in-aide no.~$14204012$. E-mail:
kitaev@pdmi.ras.ru}\\
Department of Mathematics, Kyushu University, 812-8581 Fukuoka, Japan\footnotemark[1]\\
Steklov Mathematical Institute, Fontanka 27, St. Petersburg, 191023, Russia\footnotemark[2]\\
and\\
School of Mathematics and Statistics, University of Sydney,\\
Sydney, NSW 2006, Australia\footnotemark[1]\;\,\footnotemark[2]}

\begin{document}

\maketitle

\begin{abstract}
An {almost Belyi covering} is an algebraic covering of the projective line,
such that all ramified points except one simple ramified point lie above a set of 3 points
of the projective line. In general, there are 1-dimensional families of these coverings with a fixed ramification pattern. (That is, Hurwitz spaces for these coverings are curves.)
In this paper, three almost Belyi coverings of degrees 11, 12, and 20 are explicitly constructed. 
We demonstrate how these coverings can be used for computation of several algebraic solutions
of the sixth Painlev\'e equation. 

%Additionally, we develop a rigorous approach to the theory of deformations of dessins d'enfant,
% associated to almost Belyi coverings.

%Almost Belyi coverings can be used to pullback hypergeometric equations to isomonodromic 
%$2\times2$ Fuchsian systems with 4 singularities, and to compute corresponding algebraic
%Painlev\'e VI solutions. Using the constructed three coverings, we produce five algebraic
%Painlev\'e VI solutions,  corresponding to Fucshian systems with the icosahedral monodromy group. 
%These Painlev\'e VI solutions are computed by Boalch as well, by other method.
%They have type 37, 38, 41, 42, 43 in his classification of icosahedral Painlev\'e VI functions. 
%As an upshot, we present a method of computing projective normalization (by Schlesinger 
%transformations) of the pull-backed Fuchsian systems. 

\vspace{24pt}\noindent
\vspace{12pt}\\
{\bf Short title}: {Highly ramified coverings}
\end{abstract}

\newpage

% locust   http://images.ibsys.com/2004/1129/3956033.jpg

\section{Introduction}

Recall that a {\em Belyi function} is a rational function on an algebraic
curve with at most 3 critical values. The corresponding covering of $\PP^1$
by the algebraic curve ramifies only above (at most) 3 points. By
fractional-linear transformations, the ramification locus can be chosen to be
the set $\{0,1,\infty\}\subset\PP^1$.

According to Belyi \cite{Belyi} and Grothendieck \cite{Groth}, there are
deep relations between Belyi functions and algebraic curves defined over
$\overline{\QQ}$, and dessins d'enfant.

More generally, one can consider the set of (isomorphism classes of) all
coverings on $\PP^1$ with prescribed number of ramified points and with
prescribed ramification orders above them. Such a topological configuration
space is called a {\em Hurwitz space}. If we fix a ramification pattern for
Belyi functions, we typically have a finite set of (isomorphism classes of)
Belyi functions with the prescribed ramification pattern. If we fix the {\em
hypermap} \cite{zvonkin} of 3 permutations for the monodromy group of the
covering, the Belyi map is unique.

In this article we consider coverings of $\PP^1$ which ramify only above 4
points, and such that there is only one simple ramified point in one of the 4 
fibers. We refer to these coverings as {\em almost Belyi coverings}. 
As is known, Hurwitz spaces for coverings ramified only above 4 general points have
dimension one \cite[Proposition 3.1]{zvonkin}. In fact, any algebraic curve can be
obtained as some one-dimensional Hurwitz space (with specified monodromy
permutations) \cite{diaz}.

For an almost Belyi covering of degree $n$, let us denote its ramification pattern by
$R_4(P_1|\,P_2|\,P_3)$, where $P_1,P_2,P_3$ are 3 partitions of $n$
specifying the ramification orders above three points. 
The fourth partition is assumed to be $2+1+1+\ldots+1$. 
%By {\em the extra ramification point} we refer to the simple ramification point in the fourth fiber.
%According to \cite[Proposition 2.1]{HGBAA},
%the total number of parts in $P_1$, $P_2$, $P_3$ must be equal to $n+3$.
The similar notation for a ramification
pattern for Belyi maps is $R_3(P_1|\,P_2|\,P_3)$, as in \cite{AK1}, \cite{K4}.
%\cite{HGBAA}.

The main goal of this paper is to compute generic almost Belyi coverings
$\PP^1\to\PP^1$ with the following ramification patterns:
\begin{eqnarray} \label{eq:ramif12}
& R_4\big(3+3+3+1+1+1\;|\;2+2+2+2+2+2\;|\;5+5+2\big),\\ \label{eq:ramif11} &
R_4\big(3+3+3+1+1\;|\;2+2+2+2+2+1\;|\;5+5+1\big),\\ \label{eq:ramif20}&
R_4\big(5\!+\!5\!+\!5\!+\!5\,|\,
2\!+\!2\!+\!2\!+\!2\!+\!2\!+\!2\!+\!2\!+\!2\!+\!2\!+\!2\;|\;3\!+\!3\!+\!3\!+\!3\!+\!3\!+\!2\!+\!1\!+\!1\!+\!1\big).
\end{eqnarray}
Their degree is 12, 11 and 20, respectively. 
%Notice that the $R$-notation of this paper, as in (\ref{eq:ramif12}), uses different ordering of the boxes
%than the $RS$-notation. This is done to make the ramification patterns above $z=\infty$ and $z=0$
%better visible in the main part of the paper. 

We consider the three specific coverings because of their application to the
theory of algebraic Painlev\'e VI functions. With certain almost Belyi coverings, 
one can pull-back a hypergeometric differential equation
%or more precisely, $2\times2$ matrix hypergeometric differential equations
to a parametric isomonodromic Fuchsian equation with 4 regular singular points
plus one apparent singularity. Equivalently, one may obtain isomonodromic $2\times2$ 
matrix Fuchsian systems with 4 regular singular points.
The corresponding Painlev\'e VI solutions are algebraic. Knowing suitable
almost Belyi maps, one can construct explicit examples of algebraic
Painlev\'e VI solutions \cite{HGBAA},  \cite{AK2}, \cite{Do}, and solve the corresponding 
isomonodromic Fuchsian equations explicitly in terms of hypergeometric functions.
More generally, explicit knowledge of any Hurwitz space %such highly ramified coverings
can be similarly used to solve explicitly many types of Fuchsian systems,
such as Garnier systems \cite{K2}. 

With the three almost Belyi coverings we construct, we pull-back hypergeometric equations with the icosahedral monodromy group to isomonodromic $2\times2$ Fuchsian systems with 4 regular singular points, and the same monodromy group.  In total, we compute five corresponding algebraic Painlev\'e VI solutions. They have type 37, 38, 41, 42, 43 in Boalch's classification of {\em icosahedral} Painlev\'e VI equations. Three algebraic Painlev\'e VI solutions (of type 38, 42, 41, respectively) can be constructed immediately\footnote{Our original motivation for this work was to compute a few missing examples in early versions of \cite{Bo2} of icosahedral Painlev\'e VI functions. We did our computations for type 38, 41 solutions unaware of the sixth electronic version of \cite{Bo2}. Before the next version
of \cite{Bo2} with type 42, 43 examples appeared, we had the degree 11
covering and the corresponding type 42 solution as well. Complimentary to
\cite{Bo2}, computations of type 44--45 and 47--52 examples were done independently in
\cite{Bo4} and \cite{KV1}.}
from the almost Belyi coverings with the ramification patterns
(\ref{eq:ramif12})--(\ref{eq:ramif20}). To obtain other two algebraic Painlev\'e VI solutions,
we compute properly pull-backed corresponding Fuchsian systems explicitly.
The type 41 solution is related to the Great Dodecahedron Solution of Dubrovin-Mazzocco \cite{DM}
via an Okamoto transformation.

Efficient computation of highly ramified coverings or Hurwitz spaces are important problems
in other fields as well.  Therefore these problems attract attention of researchers.
In \cite{Couv}, a method is presented to compute one-dimensional Hurwitz spaces (for almost
Belyi maps, for example) based on the degenerations when 4 ramification loci
coalesce into 3 ramified points. In \cite{GAPhurwitz}, a computer algebra
package is presented for computing genera and monodromy groups of Hurwitz
spaces or coverings.

In the next section, we present our computational method. It is basically
the same method as described in \cite[Section 3]{hhpgtr}. Compared with the
most straightforward method with undetermined coefficients, we derive
equations of smaller degree in undetermined coefficients by using properties
of the derivatives of Belyi or almost Belyi maps. We found out that very much the
same computational method was used in \cite{Hempel}, for deriving several
rather simple Belyi coverings by hand. In Sections \ref{sec:sol12} through
\ref{sec:sol20} we present our computations of the coverings with ramification
patterns (\ref{eq:ramif12})--(\ref{eq:ramif20}). Section \ref{sec:painleve} 
demonstrates application of the three coverings to computation of algebraic Painlev\'e VI solutions.
The coverings are $R$-parts of $RS$-pullback transformations of 
hypergeometric differential equations to isomonodromic $2\times 2$ Fuchsian systems.
Computation of $RS$-pullback transformations of isomonodromic systems $2\times 2$ Fuchsian systems is discussed thoroughly in \cite{KV3}.

The authors prepared a {\sf Maple 9.0} worksheet supplementing this article,
with the formulas in {\sf Maple} input format, and demonstration of key
computations. Readers may contact the authors, or search a current website
of the first author on the internet, to access the worksheet.

\section{The computational method}

Here we briefly recall the straightforward method for computation of almost Belyi coverings
from $\PP^1$ to $\PP^1$, and present an improved method that uses differentiation.
To distinguish the two projective curves, we write the coverings as $\PP^1_x\to\PP^1_z$,
where $x$ and $z$ denote the rational parameters of the projective lines
above and below, respectively. We assume that the three %(supposedly)
ramification loci indicated in the $R_4$-notation are $z=0$, $z=1$, $z=\infty$, in  this order.
We refer to the simple ramification point in the fourth ramified fiber as
{\em the extra ramification point}. 

Let $n$ denote the degree of the covering. By the Hurwitz genus formula
\cite[Corollary IV.2.4]{harts}, the number of distinct points above
$\{0,1,\infty\}\subset\PP^1_z$ must be $n+3$
for an almost Belyi covering (and $n+2$ for a Belyi map); 
see \cite[Proposition 2.1]{HGBAA} or \cite[Lemma 2.5]{V}. 

The straightforward method to compute an almost Belyi covering
with a given ramification pattern is to write an ansatz of the form
\begin{eqnarray}  \label{eq0:ans}
\varphi(x)=\frac{F}{H},\qquad \varphi(x)-1=\frac{G}{H},
\end{eqnarray}
where $F$, $G$, $H$ are general polynomials in $x$ of the factorized form
determined by the respective partition of $n$. Specifically, the polynomials
have the form $C_0\prod_{j=1}^n P_j^j$, where $C$ is a constant
(undetermined yet), and each $P_j$ is a monic general polynomial of degree
equal to the number of parts $j$ in the respective partition. Of course,
polynomials of degree zero can be effectively skipped. To avoid redundancy,
one may assume that $H$ is a monic polynomial, and may pick 3 of the
$x$-points\footnote{Strictly speaking, the $x$-points for almost Belyi coverings
are curves, or one-dimensional branches of a generic family, 
parametrized by an isomonodromy parameter $t$ or other parameter,
since the Hurwitz spaces for almost Belyi maps are one-dimensional. 
For simplicity, we ignore the dimensions introduced by  such parameters, 
and consider a one-dimensional Hurwitz space as a generic point.}
as $x=\infty$, $x=0$ and $x=1$. Expression (\ref{eq0:ans}) leads
to the polynomial identity $F=G+H$; by expanding the polynomials and comparing
the terms to the powers of $x$ one gets a set of polynomial equations.

This straightforward method was extensively used in \cite{HGBAA}, \cite{K4}
to compute Belyi maps and almost Belyi coverings of degree up to 12. 
Those coverings were applied to compute algebraic transformations
of Gauss hypergeometric functions, or
compute algebraic Painlev\'e VI functions. However, the amount of
computations with the straightforward method grows quickly for larger $n$.
In particular, the number of variables and algebraic degree of initial
equations grow linearly with $n$. The polynomial system may have many
degenerate (or {\em parasitic} \cite{kreines}) solutions, when the rational
expression in (\ref{eq0:ans}) can be simplified to a rational function of
lower degree. Computation of higher degree coverings with the
straightforward method is hardly possible even with modern computers. In
particular, our three coverings, including the degree 11 covering,
% for (\ref{eq:ramif12})--(\ref{eq:ramif20})
were too hard to compute in reasonable time with available PC's.
%personal computers.

Equations of smaller degree for the undetermined coefficients can be
obtained by considering derivatives of $\varphi(x)$. According to Couveignes
\cite{Couv}, Fricke was probably the first to use differentiation to
investigate highly ramified maps. More recently, differentiation was used
for investigation of dessins d'enfant in \cite{shabatgb}, \cite{Schneps},
\cite{zapponi} and other works. Specifically for computational purposes,
differentiation was used in \cite{Hempel} and \cite{hhpgtr} in similar ways.

A systematic procedure for computation of Belyi maps with differentiation is
formulated in \cite[Section 3]{hhpgtr}. The main trick is to consider
logarithmic derivatives of $\varphi$ and $\varphi-1$ from (\ref{eq0:ans}).
For instance, the denominator of $\varphi\,'/\varphi$ %the logarithmic derivative of $\varphi$
is the product of all factors of $F$ and $H$, to the power 1. The numerator
is the product of all factors of $G$, with the powers diminished by 1. This
gives equations of smaller degree, and easy possibilities for elimination.
Typically, the degree of equations is diminished by the number of distinct
points in a corresponding ramified fiber. Compared with computation of Belyi maps,
the only adaptation for almost Belyi coverings is that numerators of the logarithmic derivatives
have an additional degree 1 factor coming from the extra ramified point. 

In this paper, we apply the method in \cite[Section 3]{hhpgtr} for
computation of the almost Belyi coverings with the ramification type
(\ref{eq:ramif12})--(\ref{eq:ramif20}). In the following Section, we present
the computational steps specifically for ramification pattern
(\ref{eq:ramif12}) quite in detail. Having demonstrated that example, we
present our computations for the other two coverings with less notice of the
routine steps, but we concentrate rather on additional heuristic means of
solving the obtained systems of equations. In particular, in Section
\ref{sec:sol20} we use modular methods for finding coverings with ramification
pattern (\ref{eq:ramif20}). Our examples show that computational complexity
depends not only on the degree of the covering, but also on geometric
complexity of the solutions (apparently, the geometry is more complicated
when the degree is prime), or the number of irreducible components of the
solutions. Accordingly, different heuristic tricks can be useful for
different coverings.

%Computations of that much complicated coverings are rarely described. Here,
%even with simpler equations on hand, we need to make  decisions on how to
%solve them.

In the rest of this Section, we make a few comments on the fields of
definition and dimension of Hurwitz spaces. Whether we use the
straightforward method or logarithmic derivatives, the algebraic equations
for the coverings (or Hurwitz spaces) %between undetermined coefficients
are defined over $\QQ$. For Belyi functions, the solutions (up to
fractional-linear transformations) are isolated points, generally defined
over an algebraic extension of $\QQ$. The field extension may depend on the
fractional-linear normalization of fixing the points $x=\infty$, $x=0$,
$x=1$. To be certain of a minimal $\QQ$-extension, one may choose to fix the
points of $\PP^1_x$ where the ramification order is different from other
ramification orders %at other points
in the same fiber. The central question in the theory of dessins d'enfant is
how the Galois group of $\overline{\QQ}/\QQ$ acts on Belyi coverings (with
necessarily the same ramification pattern) or dessins d'enfant.

General almost Belyi coverings are parameterized by algebraic curves
\cite{zvonkin}, \cite{HGBAA}. If preferred so, one may consider them as
one-dimensional families of almost Belyi coverings. The genus of
parameterizing curves may depend on the normalization. To get a modelling
curve of minimal genus, we may strive to fix the points with ``isolated"
ramification orders as $x=\infty$, $x=0$, $x=1$, like in the zero-dimensional
case. But the equation system can be simpler if we adopt the strategy of
choosing the points with the highest ramification orders as $x=\infty$, $x=0$,
$x=1$. We demonstrate this situation in Section \ref{sec:sol11}.

%+Ad hoc

\section{The degree 12 covering}
\label{sec:sol12}

Here we compute the generic pull-back covering with the ramification type
\begin{equation} \label{eq12:ramif}
R_4\big(3+3+3+1+1+1\;|\;2+2+2+2+2+2\;|\;5+5+2\big).
\end{equation}
Relatively speaking, this is a warm-up example.

By our conventions, the three partitions of 12 specify the ramification orders
above $z=0$, $z=1$ and $z=\infty$, respectively. We choose the simple
ramified point above $z=\infty$ as $x=\infty$, and the simple ramified
point above the fourth $z$-point as $x=0$. We do not fix $x=1$, so there
will be the torus action $x\mapsto\lambda x$ on the defining equations. The
equations are expected to be weighted-homogeneous, and the Hurwitz space (if
irreducible) should be a curve in a weighted-projective space of minimal
possible genus.

We write the ansatz
\begin{equation} \label{eq12:direct}
\varphi_{12}(x)={C_0}\frac{F^3G}{H^5}, \qquad
\varphi_{12}(x)-1={C_0}\frac{P^2}{H^5},
\end{equation}
where
\begin{eqnarray*}
F \equal x^3+a_1x^2+a_2x+a_3,\\
G \equal x^3+b_1x^2+b_2x+b_3,\\
H \equal x^2+c_1x+c_2,\\
P \equal x^6+p_1x^5+p_2x^4+p_3x^3+p_4x^2+p_5x+p_6,
\end{eqnarray*}
are polynomials whose roots are the other $x$-points above
$z\in\{0,1,\infty\}$. In particular, the roots of $F$, $H$, $P$ are the
remaining ramified points. Besides, $C_0=\lim_{x\to\infty} \varphi_{12}(x)/x^2$
%is the constant $\lim_{x\to\infty} \varphi/x^2$,
is an undetermined constant yet. The straightforward method would utilize
the following consequence of (\ref{eq12:direct}):
\begin{equation} \label{naiveeq12} \textstyle
F^3G=P^2+\frac{1}{C_0}H^5. %\frac1{C_0}G^3H,
\end{equation}

Following \cite[Section 3]{hhpgtr}, we obtain simpler equations in the
coefficients of $F$, $G$, $H$, $P$ by considering the logarithmic
derivatives of $\varphi_{12}(x)$ and $\varphi_{12}(x)-1$. It is not hard to figure out the
zeroes and poles of the logarithmic derivatives:
\begin{equation} \label{eq12:logder}
\frac{\varphi_{12}'}{\varphi_{12}}=C_1\frac{x\,P}{F\,G\,H},\qquad
\frac{(\varphi_{12}-1)'}{\varphi_{12}-1}=C_2\frac{x\,F^2}{H\,P}.
\end{equation}
Here $C_1=C_2=2$ by local considerations at $x=\infty$. One may generally
notice that if $x=\infty$ is chosen above $z=\infty$ in a setting like
(\ref{eq12:direct}), the constants in the logarithmic derivative expressions
like in (\ref{eq12:logder}) are equal to the ramification order at $x=\infty$.

Comparison of the numerators in (\ref{eq12:logder}) gives the following
identities:
\begin{equation}
2xP=3F'GH+FG'H-5FGH',\qquad 2xF^2=2P'H-5PH'.
\end{equation}
The same type of expressions is derived in computations in \cite{Hempel}.
After expanding the polynomial expressions and collecting the terms to the
powers of $x$, the first identity gives the following equations:
\begin{eqnarray} \label{eqs:eq12a}
2p_1\equal 7c_1+b_1-a_1,\nonumber\\
2p_2\equal 12c_2+6b_1c_1+4a_1c_1-2a_1b_1-4a_2,\nonumber\\
%2p_3\equal 11b_3+5c_3+9b_2c_1+7b_1c_2+3a_1b_2-a_1c_2+\ldots,\\
%a_1b_1c_1-5a_2b_1-7a_2c_1-13a_3,\\
\cdots & \cdots & \cdots\\
2p_{6}\equal
6a_1b_3c_2+4a_2b_2c_2-2a_2b_3c_1+2a_3b_1c_2-4a_3b_2c_1-10a_3b_3,\nonumber\\
0\equal 3a_2b_3c_2+a_3b_2c_2-5a_3b_3c_1.\nonumber
\end{eqnarray}
The second identity gives the equations
\begin{eqnarray} \label{eqs:eq12b}
4a_1\equal 7c_1,\nonumber\\
4a_2+2a_1^2\equal 12c_2+5c_1p_1-2p_2,\nonumber\\
%8a_3+24a_1a_2+8a_1^3\equal
%11b_3+5c_3+9b_2c_1+7b_1c_2+6b_2p_1+2c_2p_1+4b_1c_1p_1+\ldots\\
%&& +b_1p_2-c_1p_2 -4p_3, %\quad \mbox{etc.}
\cdots & \cdots & \cdots \\
2a_{3}^2\equal 4c_2p_4 -3c_1p_5-10p_6,\nonumber\\
0\equal 2c_2p_5-5c_1p_6.\nonumber
\end{eqnarray}
The new equations are sufficient, since they are derived from necessary
conditions. They have smaller algebraic degree, has less degenerate
solutions, and can be solved even by brute force with {\sf Maple}'s routine
{\sf solve}. More systematically, one may use elimination or Gr\"obner basis
techniques. The system is overdetermined, but superfluous equations only
help with Gr\"obner basis computations. As mentioned, the equations are
weighted homogeneous; specifically
\begin{equation} \label{eq20:degrees}
\deg a_j=\deg b_j=\deg c_j=\deg p_j=j.
\end{equation}
The variables $p_i$ can be directly eliminated using the first set of
equations. This can be done similarly for any covering problem with a %``clean"
fiber of only simple ramified points (with the ramification order 2) plus
possibly one non-ramified point. Also notice that the second set of
equations does not contain $b_i$'s. There is a dependence between the first
two equations in (\ref{eqs:eq12a}) and the first two equations in
(\ref{eqs:eq12b}).

A straightforward way to get the result is the following. Using first
equation in (\ref{eqs:eq12b}), we eliminate $a_1$. Then all equations are
(still) linear in the $b_i$'s and $p_i$'s. We actually get 12 linearly
independent equations in these 9 variables. These variables can be
eliminated using determinants\footnote{A brute way to eliminate the 9
variables is to pick 9 (out of the 12) equations, solve them in the 9
variables, and substitute into the remaining equations. This is equivalent
to computation of $10\times10$ determinants, with polynomial entries in
$c_1,c_2,a_2,a_3$. In the particular case, this method typically gives
equations of degree 4 or 5 in $a_3$ alone.} or syzygies\footnote{Here we
mean syzygies between the 12 vectors in the rank 9 free module over
$\QQ[c_1,c_2,a_2,a_3]$, with the vector components equal to the
corresponding coefficients to the 9 variables. % in the 12 equations. %
When the syzygies are applied to the 12 equations, the 9 variables will be
eliminated. In the particular case, we can get one equation of weighted
degree 11 in $c_1,c_2,a_2,a_3$ in this way, and 3 independent equations of
degree 12.}. We get weighted homogeneous equations in $c_1,c_2,a_2,a_3$. We
just have to eliminate one more variable\footnote{Here resultants can be
used. In our computations, even if we took two equations of degree 4--5 in
$a_3$, {\sf Maple 9.5} computed a resultant with respect to $a_3$ in 15--30
seconds.} to get an equation for the Hurwitz curve. A single such equation
in 3 weighted-homogeneous variables is likely to have large degree and
superfluous factors. To avoid investigating all factors, one may compute two
or more such equations, and consider only the common factors. It turns out
that only the following factor gives non-degenerate
solutions\footnote{Degenerate solutions are those for which the polynomials
$F$, $G$, $H$ have multiple or common roots, or have the root $x=0$.
Consequently, the factors such as $a_3$, $b_3$, $c_2$ or $c_1^2-4c_2$ can be
ignored. It might be even useful to search actively for polynomials
divisible by the degeneracy factors, so that after dividing them out we
possibly get polynomials of low degree. In Section \ref{sec:sol20}, we
systematically search for polynomials divisible by two resultants defining
degeneracy for ramification pattern (\ref{eq:ramif20}).}:
\[ %begin{equation}
160a_2^2c_1^2+6912a_2c_2^2-2256a_2c_1^2c_2-188a_2c_1^4+103680c_2^3-81936c_1^2c_2^2+20328c_1^4c_2-1421c_1^6.
\] %end{equation}
This weighted-homogeneous polynomial defines a curve of genus 0. We can
normalize $c_1=1$, and parameterize as follows:
\begin{equation}
%c_1=4,\quad
c_2=\frac{(2t+1)(5t+16)}{48\,t}, \qquad a_2=
-\frac{(2t+5)(15t^2+25t+16)}{16\,t}.
\end{equation}
Going back, we consequently find parametric expressions for $a_3$ (and
immediately for $a_1$), and then for the 9 variables $b_i$'s and $p_i$'s. To
find the constant $C_0$, we can use (\ref{naiveeq12}) evaluated at any
$x\in\QQ$. %$x\not\in\{0,1,\infty\}$.

To write the generic solution more compactly, we renormalize $x\mapsto x/4$
and multiply the polynomials $F$, $G$, $H$ by some expressions in $t$.
Here is the covering:
\begin{eqnarray} \label{eq:varphi12}
\varphi_{12}(x)=C_{12}\;\frac{F_{12}^3\,G_{12}} {H_{12}^5},
\end{eqnarray}
where
\begin{eqnarray*}
C_{12}\equal-\frac{t^2(10t^2+25t+16)}{16(3t+4)^7},\\
F_{12}\equal 2tx^3+14tx^2-2(2t+5)(15t^2+25t+16)x-(2t+5)(5t+16)(t^2+10t+6),\\
G_{12}\equal 50t(10t^2+25t+16)x^3-30t(14t^3-18t^2-105t-80)x^2\\
&&-6(5t+16)(2t+1)(20t^3+35t^2+3t-16)x-(2t+1)^2(5t+16)^2(5t^2+10t+6),\\
H_{12}\equal 3tx^2+12tx+(2t+1)(5t+16).
\end{eqnarray*}
We have 
\begin{equation} \label{eq1:varphi12}
1-\varphi_{12}(x)=\frac1{16(3t+4)^7}\frac{P_{12}^2}{H_{12}^5},
\end{equation}
where
\begin{eqnarray*}
P_{12}\equal 20t^3(10t^2+25t+16)x^6-12t^3(7t^3-184t^2-490t-320)x^5\\
&& +60t^2(4068t^3+4048t+1024+200t^5+5885t^2+1386t^4)x^4\\
&& +20t^2(37627t^4+100300t^3+137092t^2+6415t^5+250t^6+92992t+24576)x^3\\
&& -60t(5t+16)(200t^7+1260t^6+3052t^5+3248t^4+1126t^3-137t^2+256t+256)x^2\\
&& -30t(5t+16)^2(100t^7+740t^6+2289t^5+3780t^4+3600t^3+2040t^2+700t+128)x\\
&& -(2t+1)(5t+16)^3(100t^7+740t^6+2289t^5+3780t^4+3600t^3+2040t^2+700t+128).
%&& -30t(5t+16)^2L_0\,x-(2t+1)(5t+16)^3L_0,\\
%L_0 \equal 100t^7+740t^6+2289t^5+3780t^4+3600t^3+2040t^2+700t+128.
\end{eqnarray*}

\section{The degree 11 covering}
\label{sec:sol11}

Here we compute the generic pull-back covering with the ramification type
\begin{equation} \label{eq11:ramif}
R_4\big(3+3+3+1+1\;|\;2+2+2+2+2+1\;|\;5+5+1\big).
\end{equation}
To get a Hurwitz curve of minimal genus, we can choose the non-ramified
point above $z=\infty$ as $x=\infty$; the non-ramified point above $z=1$ as
$x=0$; and the simple ramified point above the fourth $z$-point as $x=1$.
Accordingly, we write the ansatz
\begin{equation} \label{eq11:direct}
\varphi_{11}(x)={C_0}\,\frac{F^3\,G}{H^5}, \qquad
\varphi_{11}(x)-1={C_0}\,\frac{x\,P^2}{H^5},
\end{equation}
where $F$, $G$, $H$, $P$ are polynomials of degree 3, 2, 2 and 5,
respectively. Zeroes and poles of the logarithmic derivatives of
$\varphi_{11}(x)$ and $\varphi_{11}(x)-1$ are easy to figure out, like in
(\ref{eq12:logder}). The method gives the following identities:
\begin{eqnarray} \label{eq11:ld1}
2(x-1)P\equal 3F'GH+FG'H-5FGH',\\
\label{eq11:ld2} 2(x-1)F^2\!\equal 2xP'H-5xPH'+PH.
\end{eqnarray}
However, the resulting equations in the coefficients of $F$, $G$, $H$, $P$
are still too complicated to solve by direct elimination or Gr\"obner basis
techniques.

Equations of smaller algebraic degree can be obtained if we adopt the
strategy to normalize the points with highest ramification orders as
$x=\infty$, $x=0$, etc. (But then the Hurwitz space can have non-minimal
genus.) Accordingly, we choose the two points of ramification order 5 (above
$z=\infty$) as $x=0$ and $x=\infty$. We choose the extra ramified point
above the fourth locus as $x=1$, as just above. Then the non-ramified
points above $z=\infty$ and $z=1$ are undetermined. We denote their location
as $x=c_1$ and $x=c_2$, respectively. That gives the following ansatz:
\begin{equation} \label{eq11:altern}
\widetilde{\varphi}_{11}(x)=\widetilde{C}_0\frac{F^3\,G}{x^5\,(x-c_1)}, \qquad
\widetilde{\varphi}_{11}(x)-1=\widetilde{C}_0\frac{P^2\,(x-c_2)}{x^5\,(x-c_1)},
\end{equation}
where
\begin{eqnarray*}
F \equal x^3+a_1x^2+a_2x+a_3,\\
G \equal x^2+b_1x+b_2,\\
P \equal x^5+p_1x^4+p_2x^3+p_3x^2+p_4x+p_5,
\end{eqnarray*}
%are polynomials whose roots are the other $x$-points above $z\in\{0,1,\infty\}$.
and $C_0=\lim_{x\to\infty} \widetilde{\varphi}_{11}(x)/x^5$ is an undetermined constant. The
fractional-linear transformation on $\PP^1_x$ from (\ref{eq11:altern}) to
(\ref{eq11:direct}) is
\begin{equation} \label{eq11:frli}
x\mapsto \frac{c_1(c_2-1)\,x+c_2(1-c_1)}{(c_2-1)\,x+1-c_1}.
\end{equation}

Consideration of logarithmic derivatives of $\widetilde{\varphi}_{11}(x)$ and
$\widetilde{\varphi}_{11}(x)-1$ give the following equations:
\begin{eqnarray} \label{eq11:ld3}
5(x-1)P\equal 3x(x-c_1)F'G+x(x-c_1)FG'-(6x-5c_1)FG,\\ \label{eq11:ld4}
5(x-1)F^2\equal 2x(x-c_1)(x-c_2)P'-(5x^2-4c_1x-6c_2x+5c_1c_2)P.
\end{eqnarray}
Like in the previous example, collecting terms to the powers of $x$ gives a
system of equations in the $a_i$'s, $b_i$'s, $c_i$'s and $p_i$'s. First we
pick up the following 8 equations: the terms to the powers 5, 4, 2, 1, 0 of
$x$ in (\ref{eq11:ld3}), and the terms to the powers 6, 1, 0 in
(\ref{eq11:ld4}). Using these equations, we eliminate the $a_i$'s and
$p_i$'s. Formally, there are 2 solution components, but the one with $a_3=0$
has to be discarded as degenerate\footnote{Here is a stepwise course of
elimination. The zeroth powers to $x$ in (\ref{eq11:ld3})--(\ref{eq11:ld4}),
that is, the substitution $x=0$, gives the equations $5p_5+5c_1b_2a_3=0$ and
$5a_3^2=5c_1c_2p_5$. It is easy to eliminate $a_3,p_5$. We must ignore
solutions with $a_3=0$, so we are left with $a_3=-c_1^2c_2b_2$. Next we
consider the coefficients to the first powers to $x$ in
(\ref{eq11:ld3})--(\ref{eq11:ld4}), and eliminate $a_2,p_4$. We get, in
particular,
\[
a_2=\frac{c_1}{11}\left(5c_1b_2+12c_2b_2-3c_1c_2b_1-10c_1c_2b_2\right).
\]
Similarly, we consider the coefficients to the highest degree (5 and 6,
respectively) to $x$ in (\ref{eq11:ld3})--(\ref{eq11:ld4}), and eliminate
$a_1,p_1$. We get, %expressions in terms of $b_i$'s and $c_i$'s. In
in particular, $a_1=\left(3b_1-12c_1-5c_2+10\right)/11$.
%\[
%a_1=\frac{1}{11}\left(3b_1-12c_1-5c_2+10\right).
%\]
Having expressed the $a_i$'s and $p_1$, $p_4$, $p_5$ just in terms of
$b_i$'s and $c_i$'s, the remaining equations are linear in $p_2$, $p_3$.
These two $p_i$'s can be eliminated using discriminants or syzygies, like in
the previous Section.}. We obtain several non-homogeneous equations in the 4
variables $b_1,b_2,c_1,c_2$. To get equations for the Hurwitz curve in two
variables, we eliminate $b_1$, $b_2$ by picking up equations of minimal
degree in them and using resultants. Like in the previous Section, we can
compute several resultant polynomials (in $c_1,c_2$ only) and consider their
common divisors as candidate models for the Hurwitz curve. It turns out that
there is possibly only one component of non-degenerate solutions. It is
described by a polynomial factor of degree 15 in $c_1,c_2$.

Let $Q$ denote the degree 15 factor. It apparently defines the Hurwitz space
for the desired almost Belyi map with chosen normalization. The degree of
$Q$ in $c_1$ alone is 12; the degree in $c_2$ is just 6. According to {\sf
Maple}'s package {\sf algcurves}, $Q$ defines a curve of genus 3. Luckily,
the curve is hyperelliptic. Improvised computations produced the following
Weierstrass model:
\begin{equation}
w^2=(3t^2+3t+2)(27t^6+71t^5+130t^4+140t^3+120t^2+64t+32).
\end{equation}
The variables $c_1$ and $c_2$ can be parameterized as
\begin{eqnarray*}
c_1\equal\frac{27t^6+67t^5+116t^4+118t^3+94t^2+46t+20+w(3t^2+2t+2)}{2(t+2)(t^2+1)(2t^2+3t+3)},\\
c_2\equal % \textstyle
\frac{1107t^{12}\!+\!7641t^{11}\!+\!P_1-w(9t^3\!+\!19t^2\!+\!13t\!+\!7)
(3t^5\!+\!15t^4\!+\!15t^3\!+\!45t^2\!+\!40t\!+\!26)}
{2(t+2)(3t+1)^3(t^2+1)(2t^2+3t+3)^2(5t^2+4t+3)},
\end{eqnarray*}
where { \footnotesize
$$P_1=26055t^{10}+59035t^9+99475t^8+130463t^7
+138619t^6+121015t^5+87870t^4+51600t^3+23798t^2+7574t+1460.$$}
%\begin{equation} \nonumber
%\scriptsize .
%\end{equation}
Using these expressions and the equations in $b_1,b_2,c_1,c_2$, we
parameterize $b_1$ and $b_2$:
\begin{eqnarray*}
b_1\equal\frac{(27t^5\!-\!45t^4\!-\!190t^3\!-\!360t^2\!-\!360t\!-\!216)
\left(54t^7-P_2+w(t\!+\!7)(6t^2\!+\!3t\!+\!2)\right)}
{32(3t+1)^3(t^2+1)(2t^2+3t+3)^2(5t^2+4t+3)},\\
b_2\equal-\frac{P_3+w(t\!+\!7)(6t^2\!+\!3t\!+\!2)
(54t^7\!-\!297t^6\!-\!682t^5\!-\!1145t^4\!-\!970t^3\!-\!712t^2\!-\!304t\!-\!104)}
{32(3t+1)^3(t^2+1)^2(2t^2+3t+3)(5t^2+4t+3)^2},
\end{eqnarray*}
%$P_2=-297t^6-682t^5-1145t^4-970t^3-712t^2-304t-104$.
where $P_2=297t^6+682t^5+1145t^4+970t^3+712t^2+304t+104$, and {\footnotesize
\begin{eqnarray*}
P_3\equal2916t^{14}+11124t^{13}+191673t^{12}+764136t^{11}+1953326t^{10}+3445832t^9+4698345t^8
+5040404t^7\\
&&+4425220t^6+3147872t^5+1833840t^4+840864t^3+301376t^2+74176t+11680.
\end{eqnarray*}}
Here we can stop computations on the hyperelliptic curve. Using
fractional-linear transformation (\ref{eq11:frli}), we can express
polynomials $H$ and $G$ in (\ref{eq11:direct}):
\begin{eqnarray}
H\equal x^2-\frac{(c_1-1)(c_1+c_2)}{c_1(c_2-1)}\,x+\frac{c_2(c_1-1)^2}{c_1(c_2-1)^2},\\
G\equal
x^2-\frac{(c_1\!-\!1)(2c_1c_2\!+c_1b_1\!+c_2b_1\!+2b_2)}{(c_2-1)(c_1^2+c_1b_1+b_2)}\,x+
\frac{(c_1\!-\!1)^2(c_2^2\!+c_2b_1\!+b_2)}{(c_2\!-1)^2(c_1^2\!+b_1c_1\!+b_2)}.
\end{eqnarray}
When we write these coefficients in terms of $t,w$, the square root $w$
conveniently disappears. Hence the coefficients are just rational functions
in $t$. One may check that the algebraic relation between the coefficients
of $H$ define an irreducible curve of degree 13 (and genus 0, as
parameterized by $t$). The projective degree of the parametrization by $t$
is 13 as well, so the parametrization is minimal.

Now we can get back to the equations induced by
(\ref{eq11:ld1})--(\ref{eq11:ld2}), but knowing parametric expressions for 4
variables. Now it is straightforward to find parametric expressions for the
remaining coefficients of $F$ in (\ref{eq11:direct}). The final expression
for $\varphi_{11}$ can be simplified by the renormalization
\begin{equation} \label{eq:hint11}
x\mapsto \frac{3t^5+15t^4+15t^3+45t^2+40t+26}{(3t^2+2t+2)(5t^2+4t+3)}\,x.
\end{equation}
Accordingly, the extra ramified point is not fixed to $x=1$ in this
normalization. The final expression can be written as:
\begin{eqnarray} \label{eq:varphi11}
\varphi_{11}(x)= C_{11}\;\frac{F_{11}^3\,G_{11}} {H_{11}^5},
\end{eqnarray}
where
\begin{eqnarray*}
C_{11}\equal -\frac{(2t^2+3t+3)(3t+1)^2}{108},\\
F_{11}\equal \textstyle
x^3-\frac{4941t^6+13122t^5+19905t^4+17820t^3+10795t^2+3962t+879}{(3t+1)^2(2t^2+3t+3)}\,x^2\\
&&\textstyle+\frac{(3t+1)(432t^5+570t^4+330t^3-265t^2-340t-151)}{2t^2+3t+3}\,x+3(2t^2+3t+3)(3t+1)^4,\\
G_{11}\equal \textstyle
(3t^2\!+\!2t\!+\!2)^2x^2+\frac{27t^{10}+270t^9+945t^8+2160t^7+2745t^6+1926t^5-5t^4-1340t^3-1440t^2-720t-216}
{(2t^2+3t+3)(3t+1)^2}x\\ &&-4(3t+1)(2t^2+3t+3),\\
H_{11}\equal \textstyle
(5t^2+4t+3)x^2+\frac{135t^6+396t^5+715t^4+790t^3+610t^2+280t+82}
{2t^2+3t+3}x+(2t^2+3t+3)(3t+1)^3.
\end{eqnarray*}
Let $P_{11}$ denote the degree 5 polynomial such that 
\begin{equation}
1-\varphi_{11}(x)=\frac{1}{108(2t^2+3t+3)}\,\frac{x\,P_{11}^2}{H_{11}^5}.
\end{equation}
We have $P_{11}=(3t^2+2t+2)(2t^2+3t+3)(3t+1)x^5+\ldots$. 
%\begin{eqnarray*}
%P_{11}\equal (3t^2+2t+2)(2t^2+3t+3)(3t+1)x^5\\ &&\textstyle
%+5\,\frac{6804t^8+22734t^7+44523t^6+55962t^5+50820t^4+32678t^3+15183t^2+4506t+774}{3t+1}\,x^4 \\ &&\textstyle
%-15\frac{920727t^{14}+5361066t^{13}+16821918t^{12}+35527914t^{11}+55893222t^{10}+68760018t^9+68279994t^8+L_1}{(2t^2+3t+3)(3t+1)^3}x^3 \\ &&\textstyle
%-5\,\frac{83106t^{15}+1195560t^{14}+6891885t^{13}+24448635t^{12}+60803055t^{11}+113516721t^{10}+165097665t^9+L_2}{(2t^2+3t+3)^2}x^2 \\ && -5L_3(3t+1)^3x %\\ && 
%-9(3t^5+15t^4+15t^3+45t^2+40t+26)(2t^2+3t+3)^2(3t+1)^6,
%\end{eqnarray*}
%{\scriptsize \begin{eqnarray*}
%L_1\equal 55733262t^7+37826301t^6+21378648t^5+9992628t^4+3770224t^3+1100102t^2+224196t+25908,\\
%L_2\equal 191305215t^8+178753865t^7+135369245t^6+82903962t^5+40666090t^4+15625390t^3+4512690t^2+893520t+98388,\\
%L_3\equal 2916t^{10}+11745t^9+36045t^8+84375t^7+148905t^6+192339t^5+182565t^4+124565t^3+60375t^2+18960t+3514.
%\end{eqnarray*} }

As we see, the degree 11 covering is more complicated than the degree 12
covering of Section \ref{sec:sol12}. Apparently the geometry of prime degree
coverings is more complex.

\section{The degree 20 coverings}
\label{sec:sol20}

Here we compute generic pull-back coverings with the ramification type
\begin{equation} \label{eq20:ramif}
R_4\big(5\!+\!5\!+\!5\!+\!5\,|\,2\!+\!2\!+\!2\!+\!2\!+\!2\!+\!2\!+\!2\!+\!2\!+\!2\!+\!2
\,|\,3\!+\!3\!+\!3\!+\!3\!+\!3\!+\!2\!+\!1\!+\!1\!+\!1\big).
\end{equation}
By fractional-linear transformations, we fix the simple ramified point
above $z=\infty$ as $x=\infty$, and we choose the extra ramified point as
$x=0$.

The ansatz is
\begin{equation} \label{eq20:direct} 
\varphi_{20}(x)=C_0\frac{F^5}{G^3H}, \qquad \varphi_{20}(x)-1=C_0\frac{P^2}{G^3H}.
\end{equation}
where
\begin{eqnarray} \label{eq20:pols} 
F \equal x^4+a_1x^3+a_2x^2+a_3x+a_4,\nonumber\\
G \equal x^5+b_1x^4+b_2x^3+b_3x^2+b_4x+b_5,\nonumber\\
\label{eq:H20}H \equal x^3+c_1x^2+c_2x+c_3,\\
P \equal x^{10}+p_1x^9+p_2x^8+\ldots+p_9x+p_{10},\nonumber
\end{eqnarray}
and $C_0=\lim_{x\to\infty} \varphi(x)/x^2$ %is the constant $\lim_{x\to\infty} \varphi/x^2$,
is an undetermined constant.
%\begin{equation} \label{naiveeq}
%C_0F^5=C_0P^2+G^3H. %\frac1{C_0}G^3H,
%\end{equation}
Consideration of logarithmic derivatives of $\varphi_{20}(x)$ and
$\varphi_{20}(x)-1$ give the following identities:
%\begin{equation}
%\frac{\varphi'}{\varphi}=C_1\frac{x\,P}{F\,G\,H},\qquad
%\frac{(\varphi-1)'}{\varphi-1}=C_2\frac{x\,F^4}{G\,H\,P}.
%\end{equation}
%Here $C_1=C_2=2$, the order at $x=\infty$.
\begin{equation}
2xP=5F'GH-3FG'H-FGH',\qquad 2xF^4=2P'GH-3PG'H-PGH'.
\end{equation}
After expanding and collecting terms to the powers of $x$, we get a system
of equations in the coefficients of $F$, $G$, $H$, $P$. Since we do not fix
$x=1$, the equations are weighted homogeneous, with the same grading as
formulated in (\ref{eq20:degrees}).

There are many possibilities to eliminate variables from the equation
system. For example, one may use the first 10 equations of the first
identity to eliminate all $p_i$'s. Then we can use the first equation from
the second group to eliminate $b_1$; the subsequent equation turns out to be
void. The next 3 equations allow us to eliminate $b_3$, $b_4$ and $b_5$. But
still, there are too many variables left to solve the system by force.

Our strategy is the following. We solve the equations modulo several large
primes, isolate non-degenerate solutions, and we try to lift them to the
characteristic 0. The principle aim is to derive modular polynomial
equations which characterize only non-degenerate solutions. When lifted to
$\QQ$, those equations are expected to have low degree and rather small
coefficients, because they would contain information only about the relevant
solutions. Eventually, it turns out there are a few connected components of
non-degenerate solutions; we are able to separate them on the modular level,
so that each lifted equation system describes only one connected component.
We effectively avoid intermediate computations with huge $\QQ$-coefficients,
and consider over $\QQ$ only those equation systems which describe isolated
components of the generic solution. Actually, we are able to get just a few
new polynomial equations over $\QQ$ of low degree, but that is just enough
for a breakthrough simplification of the original system.

We use the computer algebra package {\sf Singular} \cite{Singular}, well
suited for ring-theoretic manipulation modulo large primes. For a prime
number $p$, let $\fK_p$ denote the finite field with $p$ elements. We did
computations modulo these primes:
\begin{equation} \label{eq20:primes}
p\in\{ 32003, 31991, 31981, 31973, 31963 \}.
\end{equation}
The solutions we found after considering the first 4 primes, while
computations modulo 31963 were done for checking purposes only.

We use the weighted grading as formulated in (\ref{eq20:degrees}).
To be able to discard degenerate solutions, let $Q$ denote the resultant of
$F$ and $G$ with respect to $x$, and let
$Z$ denote the resultant of $F$ and $H$. % with respect to $x$.
Their weighted degree is 20 and 12, respectively. A solution of the original
equation system is degenerate if and only if $Q=0$ or $Z=0$.

For each prime number $p$ from (\ref{eq20:primes}), we do computations in
two rings:
\begin{eqnarray*}
R_1\equal\fK_p[a_1,a_2,a_3,a_4,b_2,c_1,c_2,c_3,Z],\\
R_2\equal\fK_p[a_1,a_2,a_3,a_4,b_2,c_1,c_2,c_3,Z,Q].
\end{eqnarray*}
We assume that $b_1,b_3,b_4,b_5$ and $p_i$'s % the coefficients of $P$
are eliminated from the original system. Let $J_0$ denote the graded ideal
in $R_2$ generated by the original polynomial equations (after the
elimination) and by definitions of $Z$ and $Q$. We wish to find polynomials
in $J_0$ divisible by $Z$ or $Q$, and that we could get equations of lower
degree for non-degenerate solutions by dividing such polynomials by the
factor $Z$ or $Q$.

Let $J_1$ denote the restriction of $J_0$ onto $R_1$. For the beginning, we
compute a Gr\"obner basis for $J_1$ in $R_1$ with respect to the total %(weighted)
degree reverse lexicographic ordering with
\begin{equation} \label{ordering}
a_4\succ c_3\succ a_3\succ c_2\succ b_2\succ a_2\succ c_1\succ a_1\succ Z.
\end{equation}
With this ordering, a (weighted) homogeneous polynomial is divisible by $Z$
if and only if the leading term is divisible by $Z$. Let $G_1$ denote the
Gr\"obner basis for $J_1$. Actually, we computed $G_1$ up to bounded degree
25. %Computations were done by limiting the weighted degree (by 25).
The computations were done on a Dell laptop computer with Pentium M 1700MHz
processor on the Windows XP platform.

The first element of $G_1$ divisible by $Z$ occurs in degree 24.
(Computations up this degree take 430 seconds.) Non-degenerate solutions
should satisfy the other degree 12 factor. We have two options: either use
the new degree 12 equation immediately and recompute the Gr\"obner basis
through degrees 12 to 24, or continue computations in degree 25 in the hope
of finding more elements of $G_1$ divisible by $Z$. The second option
appears to be more acceptable since its next step takes less time (510
versus 650 seconds). Besides, the more greedy strategy of using lowest
degree new polynomials immediately leads to more frequent and lengthier
recomputations of Gr\"obner bases. In general, one may try different
tactical choices when making computations modulo first few different primes,
and then use the best options when computing modulo other primes.

%We indicate timing of our computations.

\begin{table}
\begin{center} \begin{tabular}{|c|c|c|c|c|c||c|c|} \hline
Degree & \multicolumn{7}{|c|}{Gr\"obner basis (re)computations} \\
\cline{2-8} \scriptsize (Hilbert dim) &  1st & 2nd & 3rd & 4th & 5th & 6th & 7th\\ \hline\hline %
\whilb{11}{704} & \whilb{9}{613} & \whilb{9}{613} & \whilb{9}{613} & \whilb{9}{613} & \whilb{110}{508} & \whilb{113}{496} & \whilb{24}{108} \\
\whilb{12}{1020} & \whilb{12}{848} & \whilb{13}{847} & \whilb{15}{845} & \whilb{98}{762} & \whilb{158}{587} & \whilb{143}{571} & \whilb{120}{99} \\
\whilb{13}{1432} & \whilb{16}{1128} & \whilb{19}{1124} & \whilb{21}{1120} & \whilb{104}{954} & \whilb{165}{660} & \whilb{144}{643} & \whilb{138}{105} \\
\whilb{14}{1998} & \whilb{23}{1479} & \whilb{27}{1469} & \whilb{30}{1457} & \whilb{268}{959} & \whilb{169}{738} & \whilb{149}{718} & \whilb{139}{113} \\
\whilb{15}{2724} & \whilb{32}{1877} & \whilb{40}{1853} & \whilb{47}{1829} & \whilb{328}{1022} & \whilb{178}{811} & \whilb{154}{790} & \whilb{141}{119} \\
\whilb{16}{3689} & \whilb{45}{2347} & \whilb{59}{2298} & \whilb{70}{2248} & \whilb{365}{1088} & \whilb{185}{889} & \whilb{159}{865} & \whilb{147}{127} \\
\whilb{17}{4906} & \whilb{63}{2851} & \whilb{88}{2759} & \whilb{106}{2671} & \whilb{374}{1185} & \whilb{199}{962} & \whilb{169}{937} & \whilb{\bf154}{133} \\
\whilb{18}{6486} & \whilb{85}{3414} & \whilb{123}{3255} & \whilb{157}{3099} & \whilb{418}{1251} & \whilb{209}{1040} & \whilb{180}{1012} & \whilb{\bf159}{141} \\
\whilb{19}{8448} & \whilb{123}{3980} & \whilb{178}{3720} & \whilb{229}{3470} & \whilb{433}{1348} & \whilb{229}{1113} & \whilb{197}{1084} & \whilb{\bf160}{148} \\
\whilb{20}{10943} & \whilb{166}{4575} & \whilb{237}{4174} & \whilb{327}{3776} & \whilb{480}{1414} & \whilb{250}{1191} & \whilb{221}{1160} & \\
\whilb{21}{14004} & \whilb{236}{5129} & \whilb{346}{4543} & \whilb{484}{3950} & \whilb{506}{1511} & \whilb{286}{1264} & \whilb{256}{1233} &\\
\whilb{22}{17827} & \whilb{305}{5672} & \whilb{456}{4854} & \whilb{668}{3983} & \whilb{\bf591}{1577} & \whilb{\bf353}{1342} & \whilb{317}{1312} &\\
\whilb{23}{22464} & \whilb{427}{6126} & \whilb{610}{5042} & \whilb{1021}{3828} & \whilb{\bf688}{1674} & \whilb{\bf363}{1415} & \whilb{\bf329}{1390} &  \\
\whilb{24}{28173} & \whilb{\bf535}{6539} & \whilb{\bf852}{5177} & \whilb{\bf1456}{3645} & & & \whilb{\bf346}{1478} & \\
\whilb{25}{35024} & \whilb{\bf695}{6837} &&&&&& \\ \hline\hline  % 10944, 14006, 17833, 22476, 28198, 35068
New & 1+2 & 2 & 83 & 4+97 & 9+3 & 1+2 & 7+5+1 \\ \hline %
Time (s) & 941 & 918 & 1367 & 520 & 142 & 93 & 2 \\ \hline
\end{tabular} \end{center}
\caption{Classical transformations of hyperbolic hypergeometric functions}
\label{figtabc}
\end{table}
Our computations are summarized in Table \ref{figtabc}. Recall that the
Hilbert series of a graded ring $R$ is the series $\sum_{j=0}^{\infty}
h_jt^j$, where $h_j$ is the dimension (over the ground field) of the $j$th
graded part of $R$. We refer to the numbers $h_j$ as Hilbert dimensions.

In the first column of Table \ref{figtabc}, we list the weighted degrees
from 11 to 25. The numbers in small font are the Hilbert dimensions for
$R_1$ (i.e., the number of monomials of the weighted degrees in $R_1$.) In
the second column, we give the number of elements of $G_1$ up to each
degree, and (in the small font) the Hilbert dimensions for the ring
$R_1/J_1$. As mentioned, there is one $Z$-multiple in degree 24 and two
independent $Z$-multiples on degree 25; this is indicated in the
next-to-last row, and by the bold face of the numbers above. Therefore we
have 3 new suitable equations of degree 12 and 13. Let $J_2$ denote the
ideal generated by $J_1$ and the 3 new equations.

Subsequently, we compute the Gr\"obner basis for $J_2$ up to degree 24.
Column 3 of Table \ref{figtabc} gives the same size statistics for $J_2$.
This second run gives 2 new equations of degree 12. We iterate the procedure
of adjoining new equations and recomputing Gr\"obner basis until we don't
see Gr\"obner basis elements with the leading monomial divisible by $Z$. As
indicated by Table \ref{figtabc}, the subsequent two runs give bonanzas of
83 and 101 new equations of degree 10 to 12. The same quantity and degree of
new equations occurs modulo each chosen prime, which is a good indication.
The statistics of Table \ref{figtabc} indicate the complexity of
computations in each run. In particular, the Hilbert dimensions indicate the
size of Gr\"obner basis elements in each degree.

After the 6th run, we don't get Gr\"obner basis elements divisible by $Z$.
Then we redo this run in the ring $R_2$; we use the same ordering as in
(\ref{ordering}) with additionally $Z\succ Q$. With this ordering, a
homogeneous polynomial is divisible by $Q$ if and only if the leading term
is divisible by $Q$. If the degree of the polynomial is less than $\deg
Q=20$, we still have the same criterium for divisibility by $Z$. We do the
computation in $R_2$ up to degree 24, and get one Gr\"obner basis element of
degree 23 and two elements of degree 24 divisible by $Q$. After dropping the
factor $Q$, we get an equation of degree 3 and two equations of degree 4 for
non-degenerate solutions! We feed these equations into the current Gr\"obner
basis in $R_2$; the following run gives a few equations divisible by $Z$ of
degree up to 19. After adding the new equations (of degree up to $19-12=7$),
we get a stable Gr\"obner basis: 17 equations of degree up to 7, plus
expressions for $Z$ and $Q$. This looks like the wanted polynomial equations
for non-degenerate solutions.

The computations with {\sf Singular} were semiautomatic. For the first two
primes, we tried a couple of strategies when to recompute a Gr\"obner basis.
For computations modulo other primes we chose the most effective degree
bounds. Those computations were automatic in principle; in particular, the
6th run was done directly in $R_2$. As mentioned, the new equations in
corresponding recomputations have the same degree (and leading monomials)
modulo each chosen prime. Computations for one prime take about an hour.

We would like to lift the simplest new equations to characteristic 0.
Modular lifting of Gr\"obner bases is considered by several authors; see for
example \cite{Arn} and further references. Strictly speaking, we do not know
good bounds for the size of the $\QQ$-coefficients. Sure bounds must be
huge, due to suspected complexity of intermediate computations purely in
$\QQ$. Our whole idea is to escape the intermediate computations by
considering over $\QQ$ exclusively equations for non-degenerate solutions
only. Eventually, we get actually correct solutions after using just the few
primes and a few new lifted equations. At the end of this Section we
indicate a way to check for certain that there are no other solutions.

A straightforward way to lift equations to characteristic 0 is to lift the
coefficients to rational numbers with smallest numerators and denominators.
Modular reconstruction of rational numbers is a well-known problem; the
basic algorithm is proved in \cite{wangrat}.
%We used equivalent computations of LLL reduction of the lattice generated by
A good indication that a lift of a polynomial is correct is that the
denominators of the coefficients have the same factors. If the
straightforward method did not work, we tried the LLL reduction algorithm
\cite{LLL} for the lattice generated by: the vector whose entries are the
integer coefficients modulo the composite modulus, plus one extra zero
component; and the vectors with two non-zero entries --- the entry 1 at the
extra component, and the composite modulus at some (and each) other
component. In successful cases, the shortest LLL basis vector was shorter by
several orders of magnitude than the other vectors, which is a convincing
indication. We used the extra vector component because {\sf Maple}'s {\sf
LLL} routine requires linearly independent input.

As mentioned, each 6th run of modular computations gives one equation of
(weighted) degree 3 and two equations of degree 4. The degree 3 polynomials
can be lifted quite easily. Up to a scalar factor, the result is
convincingly
\begin{eqnarray} \textstyle\nonumber
\frac{625}3c_3-\frac{875}4a_3+11c_2c_1+22b_2c_1-\frac{643}4a_2c_1
-\frac{99}{10}c_1^3-\frac{179}4a_1c_2+\frac{44}3a_1b_2\;\\
\textstyle +\frac{369}2a_2a_1
+\frac{143}{10}a_1c_1^2+\frac{3811}{60}a_1^2c_1-\frac{895}{12}a_1^3.
\end{eqnarray}
The two polynomials of degree 4 can be reliably lifted from the obtained
modular data as well. With the lifted polynomials, we can immediately
eliminate $c_3$ and $a_4$. But further brute force computations seem to be
too cumbersome still.

The final Gr\"obner basis for non-degenerate solutions appears to have new
equations of degree 5 or higher, but the modular data %from the four primes in (\ref{eq20:primes})
is not sufficient to lift them. We proceeded then to eliminate variables
modulo the primes, so to get modular equations in a minimal number of
variables. These equations indicated that there are a few irreducible
components of non-degenerate solutions. It was also reasonable to expect
that these equations would be easier to lift, even if of higher degree.

Specifically, we used Gr\"obner basis computations with respect to an
elimination ordering, still in {\sf Singular}. The Hurwitz space should be a
curve (possibly reducible), so we must have weighted homogeneous equations
in three variables. The equation in %(only) % the three variables
$a_1$, $c_1$, $a_2$ has degree 14 modulo each prime in (\ref{eq20:primes}).
It factors as follows: two factors of degree 8 and 6 modulo 31991 (and
31963); or three factors of degree 8, 3, 3 modulo other primes. This
suggests that there are three families of non-degenerate solutions, two of
them are conjugate over a quadratic extension of $\QQ$. It was possible to
lift the degree 8 and degree $6=3+3$ factors to $\QQ$ using the LLL
algorithm. %: we got two polynomial equations of degree 8 and 6 over $\QQ$.

The degree 6 polynomial factors over $\QQ(\sqrt{-15})$ into two factors of
degree 3, like expected. One may check that the factors define genus 0 %zero
curves. They can be parameterized, and the other variables can be uniquely
parameterized as well using the original equations and the three lifted
equations of degree 3 or 4. But one may notice a shortcut: ramification type
(\ref{eq20:ramif}) can be realized by a composition of two coverings of the
types:
\begin{equation}
R_3(3+\widehat{1}+\widehat{1}\,|\,5\,|\,2+2+\widehat{1}\,)
\qquad\mbox{and}\qquad R_4(3+1\,|\,2+1+1\,|\,2+2).
\end{equation}
The hats in the first expression indicate the ramification locus of the
subsequent degree 4 covering. The coverings of these two types are known due
to their application to Gauss hypergeometric functions and algebraic
Painlev\'e VI functions, as we recall in Section \ref{sec:painleve} below.
The (normalized generic) coverings for these ramification types are:
\begin{eqnarray} \label{eq:hpg35}
\varphi_5(x)\equal\frac{(5\!-\!3\sqrt{-15})\,(128x+7+33\sqrt{-15})^5}
{8000\,x\,(1024x\!-\!781\!-\!171\sqrt{-15})^3},\\ 
\label{eq:varphi4} \varphi_4(z)\equal \frac{(t-3)^3(3t-1)^3(z+1)(z+t)}
{(t\!-\!1)^2\left((t\!-\!1)^2z+t(t\!+\!1)\right)\left(4z+3(t\!+\!1)\right)^3}.
\end{eqnarray}
Note the appearance of $\sqrt{-15}$ in the first covering. The second
covering has the parameter $t$ as it is an almost Belyi function. The
composition $\varphi_4\circ\varphi_5(x)$ can be written in the form %(\ref{eq20:direct}).
(\ref{eq20:direct}). One may check that the coefficients $a_1,c_1,a_2$ in
this form parameterize the degree 6 polynomial factor (or one of the degree
3 factors). It follows that the two coverings implied by the degree 6
factor are $\varphi_4\circ\varphi_5$ and the conjugated version.

% If we leave $a_1,c_1,c_2$ we get degree 8 and 3 equation...
%-25/4*c[2]*c[1]+75/32*c[2]*a[1]-63/8*c[1]^3+109/8*c[1]^2*a[1]-27/4*c[1]*a[1]^2+a[1]^3

The degree 8 factor in $a_1$, $c_1$, $a_2$ defines a genus 0 curve as well.
Its parametrization gives rise (by the original equations and the three
lifted equations of degree 3 or 4) to the following impressive solution of
the covering problem:
\begin{eqnarray} \label{eq:coeff20}
a_1\equal 4(3t^4+9t^3+43t^2+40t+12),\\
a_2\equal 6(182t^6+728t^5+2373t^4+3584t^3+2632t^2+960t+144),\nonumber \\
a_3\equal 4(t+1)(8t^2+7t+2)(1029t^5+3246t^4+9608t^3+10224t^2+4752t+864),\nonumber\\
a_4\equal (5t+6)(8t^2+7t+2)^2(1029t^5+3246t^4+9608t^3+10224t^2+4752t+864),\nonumber\\
b_1\equal\textstyle -\frac15(4t^4+12t^3-869t^2-892t-276),\nonumber\\
b_2\equal\textstyle -\frac15(448t^6+1756t^5-59175t^4-123710t^3-101900t^2-39144t-6048),\nonumber\\
b_3\equal\textstyle -\frac15(8t^2+7t+2)\nonumber\\
&& \times(2352t^6+9399t^5-256240t^4-577360t^3-506160t^2-204336t-32832),\nonumber\\
b_4\equal\textstyle -\frac25(8t^2+7t+2)^2(56t^2+79t+34)(49t^4+132t^3-5184t^2-4752t-1296),\nonumber\\
b_5\equal\textstyle -\frac25(8t^2+7t+2)^3(49t^2+76t+36)(49t^4+132t^3-5184t^2-4752t-1296),\nonumber\\
c_1\equal\textstyle \frac13(100t^4+300t^3+507t^2+388t+108),\nonumber\\
c_2\equal\textstyle \frac13(5t+6)(8t^2+7t+2)(140t^3+307t^2+308t+108),\nonumber\\
c_3\equal\textstyle \frac13(5t+6)^2(8t^2+7t+2)^2(49t^2+76t+36).\nonumber
\end{eqnarray}
From now on, let us denote by $\varphi_{20}(x)$ the covering defined by (\ref{eq20:direct}), 
(\ref{eq20:pols}) and (\ref{eq:coeff20}). 
%\begin{equation} \label{eq:varphi20} 
%\varphi_{20}(x)=C{20}\frac{F_{20}^5}{G_{20}^3H_{20}},
%\qquad \varphi_{20}(x)-1=C_{20}\frac{P_{20}^2}{G_{20}^3H_{20}}.
%\end{equation}
We refer to the other two solutions as compositions of the degree 4 and 5 coverings.

We looked for the degree 20 coverings by modular methods, deliberately ignoring the size of intermediate would be computations over $\QQ$. It is theoretically possible that some generic
solutions are missing, since they coincide with the derived solutions modulo each of the considered primes, or the original equation system is insufficient modulo those primes. This possibility has extremely low probability. Very likely, the ``bad" primes can be only a few small prime numbers.
Regarding the application to Painlev\'e VI functions of the next section, we already know that there
are exactly two algebraic solutions of  $P_{VI}(0,0,0,-2/3;t)$ from the work of Dubrovin-Mazzocco \cite{DM}: the Cube solution and the Great Dodecahedron solution. After an Okamoto 
transformation we have exactly two solutions of $P_{VI}(1/3,1/3,1/3,1/3;t)$; they are obtainable from
the composition $\varphi_4\circ\varphi_5(x)$ and the irreducible covering $\varphi_{20}(x)$, 
as we will see.

To check for certain that there are no other covering with the ramification pattern (\ref{eq20:ramif}),
the canonical method is combinatorial. Each branch of the Hurwitz space corresponds to a 4-tuple
of permuations of 20 elements, of cycle types $5+5+5+5$, $3+3+3+3+3+2+1+1+1$, etc.
The braid group on 4 braids acts on the branches of the same connected component.
There must be only three orbits of the braid group, giving three connected components
of the Hurwitz space, corresponding to the composition $\varphi_4\circ\varphi_5(x)$, 
its complex conjugate, and $\varphi_{20}(x)$. This method is strict \cite{Couv}, \cite{GAPhurwitz}, 
but it requires computation of all permutation combinations with the given cycle type
and identity product. More geometrically, one may introduce {\em deformations of dessins d'enfant}
\cite{HGBAA}, \cite{K4}, cacti \cite{cacti} or similar geometric objects \cite[pg.~105]{Bo3} that represent
almost Belyi coverings in the same way as usual dessins d'enfant correspond to Belyi maps,
observe homotopic action of the braid group, and count possible ``deformation" drawings with the given branching type.

\section{Application to algebraic Painlev\'e VI functions}
\label{sec:painleve}

As noticed in \cite{V}, \cite{hhpgtr} and \cite{HGBAA}, \cite{AK1}, certain
Belyi coverings occur with algebraic transformations of Gauss hypergeometric
solutions. These transformations are induced by pull-back transformation of a hypergeometric differential equation to a hypergeometric equation again. 
%To keep the number of singularities of the pull-backed equation down to 3,
% the pull-back covering must be a Belyi map typically \cite[Lemma 2.2]{hhpgtr}.
In particular, Belyi covering (\ref{eq:hpg35})  transforms between
standard hypergeometric equations with the icosahedral and tetrahedral monodromy groups. 
Here is an induced hypergeometric identity:
\begin{equation}
\hpg21{1/4,-1/12}{2/3}{\,x}=\left(
1\!+\!\frac{7-33\sqrt{-15}}{128}\,x\right)^{-1/12}
\hpg21{11/60,-1/60}{2/3}{\,\frac{1}{\varphi_5(x)}}.
\end{equation}
This is the same transformation as formula (50) in \cite{V}, but with a
different definition of $\varphi_5(x)$. This formula can be checked by
comparing the Taylor expansions of both sides around $x=0$.

Similarly \cite{HGBAA}, \cite{Do}, almost Belyi coverings with certain ramification patterns can be used to pullback hypergeometric differential equations % (in a $2\times 2$ matrix form) 
to $2\times 2$ isomonodromic Fuchsian systems with four singularities. Correspondingly, 
one may derive algebraic solutions $y(T)$ of the sixth Painlev\'e equation:
\begin{eqnarray}
 \label{eq:P6}
\frac{d^2y}{dT^2}&=&\frac 12\left(\frac 1y+\frac 1{y-1}+\frac 1{y-T}\right)
\left(\frac{dy}{dT}\right)^2-\left(\frac 1T+\frac 1{T-1}+\frac 1{y-T}\right)
\frac{dy}{dT}\nonumber\\
&+&\frac{y(y-1)(y-T)}{T^2(T-1)^2}\left(\alpha+\beta\frac T{y^2}+
\gamma\frac{T-1}{(y-1)^2}+\delta\frac{T(T-1)}{(y-T)^2}\right),
\end{eqnarray}
where $\alpha,\,\beta,\,\gamma,\,\delta\in\CC$ are parameters.
The standard correspondence between solutions of the sixth Painlev\'e equation and
the mentioned isomonodromic Fuchsian systems is due to Jimbo-Miwa \cite{JM}. 
If the singular points of the Fuchsian system are $x=0$, $x=1$, $x=T$, $x=\infty$,
and the local monodromy differences at them are, respectively,
$\theta_0$, $\theta_1$, $\theta_t$, $\theta_\infty$, then the corresponding Painlev\'e equation
has the parameters
\begin{equation}
 \label{eq:para}
\alpha=\frac{(\theta_\infty-1)^2}2,\quad
\beta=-\frac{{\theta}_0^2}2,\quad\gamma=\frac{{\theta}_1^2}2,
\quad\delta=\frac{1-{\theta}_T^2}2.
\end{equation}
We denote the corresponding Painlev\'e VI equation by 
$P_{VI}(\theta_0,\theta_1,\theta_T,\theta_\infty;T)$.

General pullback transformations of $2\times 2$ Fuchsian systems $d\Psi(z)/dz=M(z)\Psi(z)$
have the following form:
\begin{equation}  \label{eq:rstrans}
z\mapsto R(x),\qquad \Psi(z)\mapsto S(x)\,\Psi(R(x)),
\end{equation}
where $R(x)$ is a rational function of $x$, and $S(x)$ is a Schlesinger transformation, usually designed to remove apparent singularities. For transformations to parametric 
isomonodromic equations, $R(x)$ and $S(x)$ may depend algebraically on parameter(s) as well. 
In \cite{K2}, \cite{HGBAA}, \cite{KV3}, these pullback transformations are called {\em $RS$-pullback transformations}, meaning that they are compositions of a rational change of the independent variable $z\mapsto R(x)$ and the Schlesinger transformation $S(x)$. The Schlesinger transformation $S(x)$ is analogous here to {\em projective equivalence} transformations $y(x)\to\theta(x)y(x)$ of ordinary differential equations.
If $S(x)$ is the identity transformation, we have a {\em direct pullback} of a Fuchsian equation.

If $z=R(x)$ is an almost Belyi covering with a suitable ramification pattern, one can pick up hypergeometric equations (in a correspondingly normalized matrix $2\times2$ form) and choose
appropriate Schlesinger transformations $S(x)$ so that the pullbacked Fuchsian equation would be isomonodromic and have four singular points, and there would be a corresponding algebraic solution of the sixth Painlev\'e equation. These $RS$-%pullback 
transformations are defined in \cite{HGBAA}, 
\cite{K4}; their algorithmic construction is considered thoroughly in \cite{KV3}.
The notation for suitable classes of these $RS$-pullback transformations is
\begin{equation}
RS^2_4\left( \,e_0\, \atop \,P_0\, \right| 
{\,e_1\, \atop \,P_1\, } \left| \,e_\infty\, \atop \,P_\infty\, \right).
\end{equation}
Here the subscripts 2 and 4 indicate a second order Fuchsian system with 4 singular points 
after the $RS$-pullback; $P_0$, $P_1$, $P_\infty$ define the ramification pattern 
$R_4(P_0\,|\,P_1\,|\,P_\infty)$ of the almost Belyi covering $R(x)$; and 
$e_0$, $e_1$, $e_\infty$ are the local exponent difference of the hypergeometric equation.

With the almost Belyi coverings $\varphi_{12}(x)$, $\varphi_{11}(x)$, $\varphi_{20}(x)$
of this paper,  we can construct $RS$-transformations of the types
\begin{eqnarray*} \textstyle
RS^2_4\left( 1/3 \atop 3+3+3+1+1+1 \right| {1/2 \atop 2+2+2+2+2+2} \left| 1/5 \atop 5+5+2 \right),\quad
RS^2_4\left( 1/3 \atop 3+3+3+1+1+1 \right| {1/2 \atop 2+2+2+2+2+2} \left| 2/5 \atop 5+5+2 \right),\\
\textstyle
RS^2_4\left( 1/3 \atop 3+3+3+1+1 \right| {1/2 \atop 2+2+2+2+2+1} \left| 1/5 \atop 5+5+1 \right),\qquad
RS^2_4\left( 1/3 \atop 3+3+3+1+1 \right| {1/2 \atop 2+2+2+2+2+1} \left| 2/5 \atop 5+5+1 \right),\\
\textstyle RS^2_4\left( 1/5 \atop 5+5+5+5 \right| {1/2 \atop 2+2+2+2+2+2+2+2+2+2} \left| 
1/3 \atop 3+3+3+3+3+2+1+1+1 \right),
\end{eqnarray*}
and derive algebraic solutions of, respectively,
\begin{eqnarray} \label{eq:p6eqs}
P_{VI}(1/3,1/3,1/3,3/5;T),  &&  P_{VI}(1/3,1/3,1/3,1/5;T), \\
P_{VI}(1/3,1/3,1/2,4/5;T),  && P_{VI}(1/3,1/3,1/2,2/5;T), \\
\label{eq:p6eqz} && P_{VI}(1/3,1/3,1/3,1/3;T).
%\\  % P_{VI}(1/3,1/3,1/3,1/5;t), &  P_{VI}(1/3,1/3,1/2,2/5;t).
\end{eqnarray}
All of the $RS$-pullbacks transform hypergeometric equations with the icosahedral
monodromy group to isomonodromic Fuchsian systems with four singular points and
the same monodromy group. The Painlev\'e VI solutions are called {\em icosahedral} \cite{Bo2};
there are 52 types of them up to branching representation of the icosahedral monodromy group, or Okamoto transformations. The solutions of (\ref{eq:p6eqs})--(\ref{eq:p6eqz}) have following Boalch types, respectively: 38, 37, 42, 43, 41. As mentioned in a footnote to the introduction section, our computations of these solutions by the method of  $RS$-transformations is independent from \cite{Bo2}.

Direct results relating $RS$-pullback transformations to algebraic Painlev\'e VI solutions are presented in \cite{KV3}. The most conveninet results are reproduced here. 
\begin{theorem} \label{kit:method}
Let $k_0,k_1,k_\infty$ denote three integers, all $\ge 2$. Let $\varphi:\PP^1_x\to\PP^1_z$ denote
an almost Belyi map, dependent on a parameter $T$. Suppose that the following conditions are satisfied:
\begin{itemize}
\item[(i)] The covering $z=\varphi(x)$ is ramified above the points $z=0$, $z=1$, $z=\infty$; there is
one simply ramified point $x=y$ above $\PP_z^1\setminus\{0,1,\infty\}$; 
and there are no other ramified points.
\item[(ii)] The points $x=0$, $x=1$, $x=\infty$, $x=T$ %are constantly 
lie above the set $\{0,1,\infty\}\subset\PP^1_z$.
 \item[(iii)] The points in $\varphi^{-1}(0)\setminus\{0,1,T,\infty\}$ are all ramified with the order $k_0$.
The points in $\varphi^{-1}(1)\setminus\{0,1,T,\infty\}$ are all ramified with the order $k_1$.
The points in $\varphi^{-1}(\infty)\setminus\{0,1,T,\infty\}$ are all ramified with the order $k_\infty$.
\end{itemize}
Let $a_0,a_1,a_T,a_{\infty}$ denote the ramification orders at $x=0,1,T,\infty$, respectively. 
Then the point $x=y$, as a function of $x=T$, is an algebraic solution of 
\begin{equation} \label{eq:p6kit}
P_{VI}\left(\frac{a_0}{k_{\varphi(0)}},\frac{a_1}{k_{\varphi(1)}},
\frac{a_T}{k_{\varphi(T)}}, 1-\frac{a_{\infty}}{k_{\varphi(\infty)}};T\right).
%\varepsilon_\varphi+\frac{(-1)^{\varepsilon_\varphi}a_{\infty}}{k_{\varphi(\infty)}};t\right),
\end{equation}
%where $\varepsilon_\varphi=1$ if the degree of $\varphi$ is even, and
%$\varepsilon_\varphi=0$ if the degree of $\varphi$ is odd.
\end{theorem}
\begin{proof} This is Theorem 3.1 in \cite{KV3}. 
\end{proof}
\begin{theorem} \label{th:llentry}
Let $z=\varphi(x)$ denote a rational covering, %an almost Belyi map,
and let $F(x)$, $G(x)$, $H(x)$ denote polynomials in $x$. 
Let $E$ denote the hypergeometric equation %$E(e_0,e_1,0,e_\infty;t;z)$
with the local exponent differences $e_0$, $e_1$, $e_\infty$ at, respectively, $z=0$, $z=1$, $z=\infty$,
Suppose that the direct pullback of $E$
with respect to $\varphi(x)$ is a Fuchsian equation with the following singularities:
\begin{itemize}
\item Four singularities are $x=0$, $x=1$, $x=\infty$ and $x=T$, with the local monodromy differences 
$d_0$, $d_1$, $d_T$, $d_\infty$, respectively. The point $x=\infty$ lies above $z=\infty$.
\item All other singularities in $\PP^1_x\setminus\{0,1,T,\infty\}$ are apparent singularities.
The apparent singularities above $z=0$ (respectively, above $z=1$, $z=\infty$) are 
the roots of $F(x)=0$ (respectively, of $G(x)=0$, $H(x)=0$).
Their local monodromy differences are equal to the multiplicities of those roots. 
%corresponding roots of $F(x)=0$.
\end{itemize}
Let $\delta$ denote a non-negative integer such that $\Delta:=\deg F+\deg G+\deg H+\delta$ is even.
Suppose that $(U_2,V_2,W_2)$ %$(U_1,V_1,W_1)$ and $(U_2,V_2,W_2)$ are syzygies 
is a syzygy between the three polynomials %in $(\ref{eq:syzgb})$,
$F$, $G$, $H$, 
%\begin{equation} \label{eq:syzgb2}
%(e_1-e_0+e_\infty)F, \quad (e_0-e_1+e_\infty)G, \quad -2e_\infty(e_0+e_1-e_\infty)H.
%\end{equation}
satisfying, if $\delta=0$,
\begin{equation} \label{eq:llsd0} \textstyle
%\deg \left( FU_1-(e_0+e_1-e_\infty)HW_1 \right)<\frac{\Delta}2, \quad
\deg U_2=\frac{\Delta}2-\deg F,\qquad \deg V_2=\frac{\Delta}2-\deg G,\qquad
\deg W_2<\frac{\Delta}2-\deg H, %\mbox{and } (\ref{eq:syzdeg0}),
\end{equation}
or, if $\delta>0$,
\begin{equation} \label{eq:llsd1} \textstyle
%\deg U_1=\frac{\Delta+\delta}2-\deg F, \qquad  
\deg U_2<\frac{\Delta+\delta}2-\deg F, \quad
\deg V_2<\frac{\Delta+\delta}2-\deg G,\quad
\deg W_2=\frac{\Delta-\delta}2-\deg H.
%\mbox{and $(\ref{eq:syzdeg1})$--$(\ref{eq:syzdeg2})$},
\end{equation}
%For a function $f$, let $\LD(f)$ denote the logarithmic derivative $f'/f$.
Then the numerator of the (simplified) rational function 
\begin{eqnarray} \label{eq:llentry}
\frac{U_2W_2}{G}\!\left(\! \frac{(e_0-e_1+e_\infty)}2\frac{\varphi'}{\varphi}
-\frac{(FU_2)'}{FU_2}+\frac{(HW_2)'}{HW_2} \right) 
+\frac{(e_0-e_1-e_\infty)}2\frac{V_2W_2}{F}\frac{\varphi'}{\varphi-1} \nonumber\\
+\frac{(e_0+e_1-e_\infty)}2\frac{U_2V_2}{H}\frac{\varphi'}{\varphi(\varphi-1)},
\end{eqnarray}
has degree $1$ in $x$, and the $x$-root of it %the numerator 
is an algebraic solution of  $P_{VI}(d_0,d_1,d_T,d_\infty+\delta;t)$.
\end{theorem}
\begin{proof} This is Theorem 6.1 in \cite{KV3}. 
\end{proof}

The first theorem here is actually a special case of the second one, when the syzygy $(U_2,V_2,W_2)$ has one of the components equal to zero. The Painlev\'e VI solutions obtained in this case can be seen as inverse functions of particular projections \cite[Figure 5]{Couv} of the respective Hurwitz spaces to 
$\PP^1_t$. The implied $RS$-transformation is 
$RS^2_4\left( \,1/k_0\, \atop P_0 \right| {\,1/k_1\, \atop P_1 } \left| \,1/k_\infty\, \atop P_\infty \right)$,
where $P_0$, $P_1$, $P_\infty$ define the ramification pattern of $\varphi(x)$. The Painlev\'e VI solutions can be derived without computing the Schlesinger part $S(x)$ of the $RS$-transformation. It is these solutions that are implied or computed in  \cite{Do} and \cite{HGBAA}.

The more general Theorem \ref{th:llentry} means that the same almost Belyi covering can be used to 
pullback several hypergeometric equations, and hence derive several algebraic Painlev\'e VI equations.
For each implied $RS$-transformation, the Schlesinger transformations depend actually on two syzygies
for the same polynomial triple $(F,G,H)$, but a single Painlev\'e solution depends on one syzygy, as stated in the theorem. (The other syzygy can be usually used to compute an algebraic solution of other Painlev\'e VI equation.) 

Both theorems have to be applied to fractional-linear normalizations of the coverings $\varphi_{12}(x)$,
$\varphi_{11}(x)$, $\varphi_{20}(x)$, where three (of the four) singular points of the transformed Fuchsian equation are chosen to be $x=0$, $x=1$, $x=\infty$. Theorem \ref{kit:method} eventually  gives solutions of
%respectively,
%\begin{eqnarray} \label{eq:p6eqs}
$P_{VI}(1/3,1/3,1/3,3/5;T)$, % \nonumber\\
$P_{VI}(1/3,1/3,1/2,4/5;T)$, $P_{VI}(1/3,1/3,1/3,1/3;T)$. 
%\\  % P_{VI}(1/3,1/3,1/3,1/5;t), &  P_{VI}(1/3,1/3,1/2,2/5;t).
%\end{eqnarray}
These are icosahedral solutions of Boalch types 38, 42, 41,  respectively.

In particular, an $RS$-pullback 
$RS^2_4\left( 1/3 \atop 3+3+3+1+1+1 \right|  {1/2 \atop 2+2+2+2+2+2} \left| 1/5 \atop 5+5+2 \right)$
with respect to $z=\varphi_{12}(x)$ is a Fuchsian system with singularities at $x=\infty$ and  the roots of
$G_{12}(x)=0$. The local monodromy differences are, respectively,  $3/5,1/3,1/3,1/3$.
A suitable normalizing fractional-linear transformation may leave
$x=\infty$ invariant, but it must move two roots of $G_{12}(x)$ to the locations $x=0$ and $x=1$.
Finding one root of $G_{12}$ is equivalent to considering $G_{12}(x,t)=0$ as an equation
for an algebraic curve. The curve has genus 0; it can be parametrized as follows:
\begin{equation} \label{eq:par12tx}
 t=-\frac{3(2s+1)}{2(s^3+3s+1)}, \qquad 
 x=-\frac{(8s^2-2s+17)(2s^2+2s+3)(s+1)^2}{10(s^3+3s+1)(2s^2+s+2)}.
\end{equation}
After the repameterization of $t$ by $s$, the other two roots of $G_{12}(x)$ are equal to
\begin{eqnarray*}
& x=%x_{\pm}=
\frac{(s-2)(4s+1)(32s^6-8s^5+164s^4-94s^3+91s^2+2s+18)}{20(s^3+3s+1)(4s^2-s+1)^2}
\pm\frac{(s-2)^2(4s+1)^2(8s^3+6s-1)\,w}{20(2s+1)(s^3+3s+1)(4s^2-s+1)^2},
\end{eqnarray*}
where $w=\sqrt{(s-2)(2s+1)(2s^2+s+2)}$. 
It appears that the polynomial $G_{12}(x)$ splits over the function field of a genus 1 curve. 
Let us denote the $x$-root in (\ref{eq:par12tx}) by $c_0$,
and the latter 2 roots by $c_+$, $c_-$. A normalizing projective coordinate for $\PP^1_x$ is:  
%fractional-linear transformation 
\begin{equation} \label{eq:lambda}
%\tilde\varphi_{k}(\lambda)=\varphi_{k}(x),\qquad k=12,\qquad
\lambda_{12}(x)=\frac{x-c_+}{c_--c_+}.
\end{equation}
Equivalently, the normalizing fractional-linear substitution to the new coordinate is given by
$\lambda_{12}^{-1}(x)$:  $x\mapsto c_+(1-x)+c_-x$.  %$x\mapsto c_+(1-\lambda)+c_-\lambda$.
By Theorem \ref{kit:method} applied to $\varphi_{12}(\lambda_{12}^{-1}(x))$, 
an algebraic solution $y_{38}(T_{38})$ of $P_{VI}(1/3,1/3,1/3,3/5;T_{38})$ is given by 
\begin{equation} \label{eq:y38lambda}
T_{38}=\lambda_{12}(c_0),\qquad y_{38}=\lambda_{12}(0).
\end{equation}
Explicitly, we have
\begin{align*}
&T_{38}=\frac12+\frac{3(2s+1)(32s^7+32s^6+138s^5+25s^4+130s^3+30s^2+20s-10)}
{2(s-2)^2(4s+1)^2(2s^2+s+2)\sqrt{(s-2)(2s+1)(2s^2+s+2)}},\\
&y_{38}=\frac12+\frac{(2s+1)(32s^6-8s^5+164s^4-94s^3+91s^2+2s+18)}
{2(s-2)(4s+1)(8s^3+6s-1)\sqrt{(s-2)(2s+1)(2s^2+s+2)}}.
\end{align*}
To get the parametrization obtained in \cite{Bo2}, one has to change
$s\to-(s+1)/2s$ and the branch of the square root.

Similarly, to apply Theorem \ref{kit:method} to the degree 11 covering $\varphi_{11}(x)$,
we have to compose it with a fractional-linear transformation of $\PP_x^1$
which leaves $x=\infty$ invariant, and moves the roots of $G_{11}(x)$ from (\ref{eq:varphi12})
to the locations $x=0$ and $x=1$. The roots of $G_{11}(x)$ are:
\begin{eqnarray*}
 \label{eq:cpm11}
x\equal-\frac{27t^{10}\!+\!270t^9\!+\!945t^8\!+\!2160t^7\!+\!2745t^6
\!+\!1926t^5\!-\!5t^4\!-\!1340t^3\!-\!1440t^2\!-\!720t\!-\!216}
%\pm3\sqrt{3}\,w}
{2(3t+1)^2(3t^2+2t+2)^2(2t^2+3t+3)}\\
&& \pm \frac{3(t+2)^2(t^2+1)^2(3t^2+2t+2)\sqrt{3(t+2)(t+7)(3t^2+3t+2)}}
{2(3t+1)^2(3t^2+2t+2)^2(2t^2+3t+3)}.
%&w=(t+2)^2(t^2+1)^2(3t^2+3t+2)\sqrt{(3t^2+3t+2)(t+2)(t+7)}.&
% \nonumber
\end{eqnarray*}
Let us denote these two roots by $c_+$ and $c_-$. Then we can use the same expression
(\ref{eq:lambda}) for the normalizing fractional-linear transformation. To distinguish, we denote this fractional-linear transformation by $\lambda_{11}(x)$.
Now, let $y_0$ denote the extra ramification point of $\varphi_{11}$, outside the fiber of
$\{0,1,\infty\}\subset\PP_z^1$. %Recalling (\ref{eq:hint11}), 
We have: 
$$
y_0=\frac{3t^5+15t^4+15t^3+45t^2+40t+26}{(3t^2+2t+2)(5t^2+4t+3)}.
$$
By Theorem \ref{kit:method} applied to $\varphi_{11}(\lambda_{11}^{-1}(x))$,
an algebraic solution %$y_{42}(T_{42})$ 
of $P_{VI}(1/3,1/3,1/2,4/5;T_{42})$ is given by
$$
T_{42}=\lambda_{11}(0),\qquad y_{42}=\lambda_{11}(y_0),
$$
We arrive at the following Painlev\'e VI solution:
\begin{eqnarray*}
T_{42}\equal\frac12+\frac{27t^{10}\!+\!270t^9\!+\!945t^8\!+\!2160t^7\!+\!2745t^6
\!+\!1926t^5\!-\!5t^4\!-\!1340t^3\!-\!1440t^2\!-\!720t\!-\!216}
{6(t+2)^2(t^2+1)^2(3t^2+3t+2)\sqrt{3(t+2)(t+7)(3t^2+3t+2)}},\\
y_{42}\equal\frac12+\frac{(t+7)(45t^6+144t^5+258t^4+228t^3+121t^2+24t-12)}
{6(5t^2+4t+3)(t^2+1)(t+2)\sqrt{3(t+2)(t+7)(3t^2+3t+2)}}.
\end{eqnarray*}
%\begin{eqnarray*}
%&t_{42}=\frac12+\frac{27s^{10}+270s^9+945s^8+2160s^7+2745s^6+1926s^5-5s^4-1340s^3-1440s^2-720s
%-216}{6\sqrt{3}(s+2)^2(s^2+1)^2(3s^2+3s+2)\sqrt{(3s^2+3s+2)(s+2)(s+7)}},&\\
%&y_{42}=\frac12+\frac{(s+7)(45s^6+144s^5+258s^4+228s^3+121s^2+24s-12)}
%{6\sqrt{3}(5s^2+4s+3)(s^2+1)(s+2)\sqrt{(3s^2+3s+2)(s+2)(s+7)}}.\phantom{AAAAAAAAAAaa}&
%\end{eqnarray*}
The solution has genus 1 as well. To get the parametrization of the same solution
obtained in \cite{Bo2}, one may substitute $t\mapsto-(2s-1)/(s+2)$.

To get a solution of $P_{VI}(1/3,1/3,1/3,1/3; T)$ by Theorem \ref{kit:method},
the covering $\varphi_{20}(x)$ can be composed with a fractional-linear
transformation of $\PP_x^1$ which leaves $x=\infty$ invariant, and moves two roots of
$H_{20}(x)$ to the locations $x=0$ and $x=1$. 
Like in the case with $\varphi_{12}(x)$, one root of $H_{20}(x)$ can be made explicit by 
parametrizing the curve $H_{20}(x,t)=0$:
\begin{equation} \label{eq:par20tx}
 t=-\frac{2(2s^3+4s^2-4s+3)}{5(2s-1)^2}, \qquad 
 x=-\frac{16(s-2)^2(2s^2+s+2)^2(s^2-2s+6)}{75(2s-1)^4}.
\end{equation}
After the repameterization of $t$ by $s$, the other two roots of $H_{20}(x)$ are equal to
\begin{eqnarray*}
%& x=-\frac{16s(8s^2-11s+8)\left(56s^7-166s^6+318s^5-269s^4+31s^3+75s^2-28s+8
%\pm9(s-1)(2s^3+4s^2-4s+3)\sqrt{s(8s^2-11s+8)}\right)}{25(2s-1)^8},
x\equal-\frac{16s(8s^2-11s+8)(56s^7-166s^6+318s^5-269s^4+31s^3+75s^2-28s+8)}{25(2s-1)^8}\\
&& \pm \frac{144s(s-1)(8s^2-11s+8)(2s^3+4s^2-4s+3)\sqrt{s(8s^2-11s+8)}}{25(2s-1)^8},
\end{eqnarray*}
Let us denote the $x$-root in (\ref{eq:par20tx}) by $c_0$,
and the latter 2 roots by $c_+$, $c_-$. Then a suitable projective parameter $\lambda_{20}(x)$ 
is given by the same expression as on the right-hand side of (\ref{eq:lambda}).
An algebraic solution $y_{41}(T_{41})$ of  $P_{VI}(1/3,1/3,1/3,1/3;T_{41})$ is given by 
\begin{equation}
T_{41}=\lambda_{20}(c_0),\qquad y_{41}=\lambda_{20}(0).
\end{equation}
Explicitly, we have
\begin{align*}
&T_{41}=\frac12+\frac{(s\!+\!1)
(32s^8\!-\!320s^7\!+\!1112s^6\!-\!2420s^5\!+\!3167s^4\!-\!2420s^3\!+\!1112s^2\!-\!320s\!+\!32)}
{54s(s-1)\big(\sqrt{s(8s^2-11s+8)}\big)^3},\\
&y_{41}=\frac12-\frac{56s^7-166s^6+318s^5-269s^4+31s^3+75s^2-28s+8}
{18s(s-1)(2s^3+4s^2-4s+3)\sqrt{s(8s^2-11s+8)}}.
\end{align*}
To get the parametrization presented in \cite[Theorem C]{Bo2}, one has to change
$s\to1/s$ and the branch of the square root. This solution is related via an Okamoto transformation
to Great Dodecahedron Solution \cite[pages 134--143]{DM}. The Dubrovin-Mazzocco solution solves
$P_{VI}(0,0,0,-2/3;T_{41})$. % Dubrovin-Mazzocco \cite{DM}.

The two composite coverings $\varphi_4\circ\varphi_5$ of degree 20, with the same
ramification pattern as $\varphi_{20}$, generate the same algebraic Painlev\'e VI solution
as the degree 4 covering $\varphi_4$ defined in (\ref{eq:varphi4}). That algebraic solution
of the same equation $P_{VI}(1/3,1/3,1/3,1/3;T_{41})$ is parametrized in \cite[Section 3.2]{HGBAA},
%as a result of the $RS$-pullback $RS_4^2\left( 1/3 \atop 2+1+1 \right| 
% {1/3 \atop 3+1} \left| 1/2 \atop 2+2 \right)$,
basically using Theorem \ref{kit:method}. Via the same Okamoto transformation, 
we get the Cube solution in \cite{DM}, of $P_{VI}(0,0,0,-2/3;T_{41})$ as well. Composition
of an almost Belyi covering with a Belyi covering never changes the algebraic Painlev\'e VI
solution. Complimenting computations in \cite[Sections 3.1, 3.4]{HGBAA}, up to Okamoto transformations
we have all 5 Dubrovin-Mazzocco solutions in \cite{DM} now derived via the method
of $RS$-pullback transformations.

Theorem \ref{th:llentry} is needed to get solutions of the equations $P_{VI}\left(1/3,1/3,1/3,1/5;T\right)$ and
$P_{VI}\left(1/3,1/3,1/2,2/5;T\right)$ of Boalch types 37 and 43. 
The implied $RS$-transformations  are
$RS^2_4\left( 1/3 \atop 3+3+3+1+1+1 \right| {1/2 \atop 2+2+2+2+2+2} \left| 2/5 \atop 5+5+2 \right)$
and $RS^2_4\left( 1/3 \atop 3+3+3+1+1 \right| {1/2 \atop 2+2+2+2+2+1} \left| 2/5 \atop 5+5+1 \right)$.
At the end, the same normalized coverings $\varphi_{12}(\lambda_{12}^{-1}(x))$ and 
$\varphi_{11}(\lambda_{11}^{-1}(x))$ can be used.  But for intermediate computations of syzygies and application of formula (\ref{eq:llentry}), we may work with the simpler parametrized coverings 
$\varphi_{12}(x)$ and $\varphi_{11}(x)$. One convenient circumstance is that the point $x=\infty$
does not have to be moved.

In particular, the direct pullback of a hypergeometric equation with the local exponent differences $1/3,1/2,2/5$ 
with respect to the covering $z=\varphi_{12}(x)$ is a Fuchsian system with actual singularities at $x=\infty$
and the roots of $G_{12}(x)$, and apparent singularities at the roots of $F_{12}$, $P_{12}$
and $H_{12}$. The local monodromy differences at the actual singularities are $4/5$ or $1/3$, while 
those differences at the apparent singularities are equal to 1 or (at the roots of $H_{12}$) to 2. 
To get rid of apparent singularities after the (implied) Schlesinger transformation $S(x)$,
we have to compute syzygies between the polynomials $F_{12}$, $P_{12}$, $H_{12}^2$. 
To compute just a solution of $P_{VI}\left(1/3,1/3,1/3,1/5;T\right)$, we apply 
Theorem \ref{th:llentry} with $\delta=1$, $\Delta=12$.  The suitable syzygy is unique up to constant 
(in $x$) multiples: 
\begin{eqnarray}   \label{eq:syzs37} 
\left( \,t\,L_1,\, 1,\, -4(3t+4)^2 \left(2tx^2-t(t-8)x-(5t+16)(2t^2+3t+2)\right) \right),
\end{eqnarray}
where {\footnotesize
\begin{eqnarray*}
L_1\equal 2t(112t^2+307t+208)x^3-2t(-1421t+60t^3-392t^2-1040)x^2 \\
&&-2(5t+16)(90t^4+251t^3+183t^2+27t+16)x-(2t+1)(5t+16)^2(10t^3+29t^2+22t+2).
\end{eqnarray*}}
The $x$-root of expression (\ref{eq:llentry}) becomes
\begin{equation}
x_{37}=-\frac{(2t+1)(5t+4)(5t+16)}{2(11t^2-4t-16)}.
\end{equation}
After the normalization by $\lambda_{12}^{-1}$ we conclude that a solution of
$P_{VI}(1/3,1/3,1/3,1/5;T_{37})$ is parametrized by $T_{37}=\lambda_{12}(c_0)$
and $y_{37}=\lambda_{12}(x_{37})$, like the solution in (\ref{eq:lambda}).
We have the solution $y_{37}(T_{37})$ with $T_{37}=T_{38}$ and
\[ %begin{equation}
y_{37}=\frac12-\frac{(2s\!+\!1)
(256s^8\!-\!832s^7\!-\!800s^6\!-\!3232s^5\!-\!1844s^4\!-\!2950s^3\!-\!1436s^2\!-\!391s\!+\!64)}
{2w(4s+1)(64s^6+336s^4+104s^3+36s^2-132s-59)}.
\] %end{equation}
To get the parametrization obtained in \cite{Bo2}, one has to substitute
$s\to-(s+1)/2s$.

Similarly, the direct pullback of the same hypergeometric equation with respect to the covering 
$z=\varphi_{11}(x)$ has  apparent singularities at the roots of $F_{11}$, $P_{11}$
and $H_{11}$. For a solution of $P_{VI}\left(1/3,1/3,1/2,2/5;T\right)$, we have to compute syzygies
between the polynomials $F_{11}$, $P_{11}$, $H_{11}^2$, assuming $\delta=0$, $\Delta=10$.
A suitable syzygy is:
\begin{eqnarray} \label{eq:syzs43} \textstyle 
\!\left( L_2,\, x-(3t+1)^2,\, -30\frac{3t^2+3t+2}{3t+1}\left((18t^2+13t+9)x-(t-1)(3t+1)^3\right) \right),
 \end{eqnarray}
where $L_2$ can be computed knowing the other two components.
The $x$-root of expression (\ref{eq:llentry}) becomes
\begin{equation}
x_{43}=\frac{(3t+1)^2(13t^5+65t^4+165t^3+195t^2+140t+46)}
{3(3t^2+2t+2)(5t^6+30t^5+45t^4+22t^3-13t^2-16t-9)}.
\end{equation}
After the normalization by $\lambda_{11}^{-1}$ we conclude that a solution of
$P_{VI}(1/3,1/3,1/2,2/5;T_{43})$ is parametrized by
$T_{43}=\lambda_{11}(0)$, $y_{43}=\lambda_{11}(x_{43})$.  
We have the solution $y_{43}(T_{43})$ with $T_{43}=T_{42}$ and
\begin{eqnarray*}
y_{43}\equal\frac12+\frac{t+7}{18(t^2+1)}\times\\
&&\frac{135t^9+540t^8+1530t^7+2916t^6+3714t^5+3486t^4+2278t^3+1144t^2+399t+114}
{(5t^6+30t^5+45t^4+22t^3-13t^2-16t-9)\sqrt{3(t+2)(t+7)(3t^2+3t+2)}}.
\end{eqnarray*}
%with $K_6=135t^9+540t^8+1530t^7+2916t^6+3714t^5+3486t^4+2278t^3+1144t^2+399t+114$.
To get the parametrization obtained in \cite{Bo2}, one has to substitute
$t\to-(2s-1)/(s+2)$.

%In principle, the solutions $y_{37}(T_{37})$ and $y_{43}(T_{43})$ can be obtained by 
%Theorem \ref{kit:method},
%but the degree of suitable almost Belyi coverings is 24 and 17, respectively. The ramification patterns are
%$R_4(3+3+3+3+3+3+3+1+1+1\,|\,2+2+2+\ldots+2\,|\,5+5+5+5+4)$ and 
%$R_4(3+3+3+3+3+1+1\,|\,2+2+2+2+2+2+2+2+1\,|\,5+5+2)$, respectively.

%Rational representation of Modular numbers, Peter Hintenaus, Vilmar Trevisan
%On Integer-to-Rational Conversion Algorithm, T. Sasaki, M. Sasaki
%Maximal Quotient Rational Reconstruction: An Almost Optimal algorithm for Rational Reconstruction,
 %%  Michael Monagan

\end{document}